\documentclass[mathpazo]{cicp}
\input{preamble.sty}

\begin{document}
\title{FEALPy: A Cross-platform Intelligent Numerical Simulation Engine}


\author[Wei H Y et.~al.]{
    Yangyang Zheng\affil{1},
    Huayi Wei\affil{1, 3}\comma\corrauth,
    Yunqing Huang\affil{1, 2},
    Chunyu Chen\affil{4},
    Tian Tian\affil{5},
    Hanbin Liu\affil{1},
    Wenbin Wang\affil{1} and
    Liang He\affil{1}
}
\address{
    \affilnum{1}\ School of Mathematics and Computational Sciences,
    Xiangtan University,
    Hunan 411105, P.R. China
    \affilnum{2}\ National Center for Applied Mathematics in Hunan,
    Xiangtan 411105, P. R. China
    \affilnum{3}\ Hunan Key Laboratory for Computation and Simulation in Science
    and Engineering Xiangtan 411105, China
    \affilnum{4}\ School of Mathematical Sciences,
    Peking University,
    Beijing 100871, China
    \affilnum{5}\ Shenzhen International Center For Industrial And Applied Mathematics;
    Shenzhen Research Institute of Big Data;
    School of Science and Engineering, Chinese University of Hong Kong (Shenzhen);
    Guang Dong 518172, China
}
\emails{
    {\tt 202531510218@smail.xtu.edu.cn} (Y.~Zheng),
    {\tt weihuayi@xtu.edu.cn} (H.~Wei)
    {\tt huangyq@xtu.edu.cn} (Y.~Huang),
    {\tt cbtxs@math.pku.edu.cn} (C.~Chen),
    {\tt tiantian00@cuhk.edu.cn} (T.~Tian),
    {\tt 202431510132@smail.xtu.edu.cn} (H.~Liu),
    {\tt 202531510193@smail.xtu.edu.cn} (W.~Wang),
    {\tt lianghe@smail.xtu.edu.cn} (L.~He),
}

\begin{abstract}
In resent years, the software ecosystem for numerical simulation still remains
fragmented, with different algorithms and discretization methods often
implemented in isolation, each with distinct data structures and programming
conventions.
This fragmentation is compounded by the growing divide between packages from
different research fields and the lack of a unified, universal data structure,
hindering the development of integrated, cross-platform solutions.
In this work, we introduce FEALPy, a numerical simulation engine built around a
unified tensor abstraction layer in a modular design.
It enables seamless integration between diverse numerical methods along with
deep learning workflows.
By supporting multiple computational backends such as NumPy, PyTorch, and JAX,
FEALPy ensures consistent adaptability across CPU and GPU hardware systems.
Its modular architecture facilitates the entire simulation pipeline, from mesh
handling and assembly to solver execution, with built-in support for automatic
differentiation.
In this paper, the versatility and efficacy of the framework are demonstrated
through applications spanning linear elasticity, high-order PDEs, moving mesh
methods, inverse problems and path planning.
\end{abstract}


\maketitle

\section{Introduction}
\label{sec:introduction}

Numerical computing software plays a central role in solving a wide range of
scientific and engineering problems.
Many such problems are governed by partial differential equations (PDEs), whose
solution spaces are inherently infinite-dimensional.
To make these problems computable on digital hardware with finite memory,
numerical discretization methods—such as the finite element method (FEM),
finite difference method (FDM), and finite volume method (FVM)—approximate
infinite-dimensional spaces with finite ones.
The success of these methods has enabled accurate simulations of physical
systems in diverse fields, from fluid dynamics to solid mechanics.

However, the software ecosystem for numerical simulation remains fragmented.
Different algorithms are often implemented independently, each with its own
data structures, programming standards, and interface conventions.
At the same time, the growing diversity of computing hardware—from CPUs and GPUs
to specialized accelerators—introduces additional complexity.
Algorithms must be adapted to distinct programming models and memory
hierarchies, making cross-platform compatibility and performance portability
major challenges.
As a result, reproducing and extending numerical methods across heterogeneous
environments is cumbersome and error-prone.

In recent years, deep learning (DL) has brought new paradigms to scientific
computing.
Neural networks (NNs) can approximate complex mappings from data, offering a
data-driven complement to traditional physics-based approaches.
Yet, traditional numerical software and DL frameworks remain largely
disconnected—most numerical solvers rely on mesh-based data structures, while
DL systems such as PyTorch and TensorFlow operate primarily on tensors.
This incompatibility hinders seamless integration between the two paradigms.

Interestingly, many core operations in numerical methods can naturally be
expressed as tensor computations.
Quantities such as fields, basis functions, and element matrices can all be
represented as multi-dimensional arrays attached to mesh entities.
Fundamental operations—indexing, reduction, concatenation, and Einstein
summation—are shared between numerical algorithms and neural network training.
This theoretical overlap opens a path toward unifying traditional numerical
computation and modern tensor-based DL frameworks.

Mature numerical PDE software such as iFEM\cite{Chen:2008ifem},
AFEPack\cite{Cai:2024afepack}, MFEM\cite{Anderson:2021mfem},
FEniCSx\cite{BarattaEtal2023, ScroggsEtal2022, BasixJoss, AlnaesEtal2014},
deal.II\cite{Arndt:2019dealII}, PHG,
and OpenFOAM\cite{Jasak:2009OpenFOAM}
offer powerful and well-established infrastructures for numerical simulation.
These software systems have greatly advanced the development and deployment of
high-performance PDE solvers across many scientific and engineering domains.
While emerging demands—such as unified tensor representations, interaction with
modern DL toolchains, and cross-backend portability—call for complementary
approaches that extend beyond the traditional designs of existing packages.

Motivated by this observation, we have developed FEALPy, a tensor-centric
numerical simulation engine designed for flexibility, interoperability, and
extensibility.
FEALPy builds a unified abstraction layer over multiple tensor computation
backends (including NumPy\cite{HARRIS2020-NUMPY},
PyTorch\cite{PASZKE2019-PYTORCH}, JAX\cite{JAX2018} and more in the future),
enabling consistent data flow across heterogeneous hardware architectures.
Its modular design decouples low-level tensor operations from higher-level
algorithmic modules—such as meshes, function spaces, discretization schemes,
and solvers—allowing users to implement complex simulation workflows with
concise and reusable code.
By adopting a unified tensor interface, FEALPy provides a common computational
foundation for both classical numerical algorithms and differentiable
programming.
It supports automatic differentiation (AD) throughout the entire pipeline,
making it possible to embed finite element solvers as differentiable layers
within NNs.
This integration bridges the gap between traditional physics-based modeling and
data-driven learning, advancing the convergence of numerical simulation and
modern AI.

FEALPy is fully open-source under the GNU GPL v3.0 license, with code hosted on
GitHub at \url{https://github.com/weihuayi/fealpy}.

The remainder of this paper is organized as follows.
Section 2 introduces the overall software architecture and data structure of
FEALPy.
Section 3 discusses implementation details and programming techniques.
Section 4 demonstrates applications of FEALPy in diverse fields, including
elasticity, high-order PDEs, moving mesh, operator learning, and optimization.
Section 5 presents two application frameworks, SOPTX and FractureX, developed
on top of FEALPy.
Section 6 concludes the paper by summarizing the main contributions and
discussing the future directions.

\section{Software Architecture and Data Structure}
\label{sec:software_architecture}

\subsection{Overview}

\begin{figure}[htbp]
    \centering
    \includegraphics[width=0.95\textwidth]{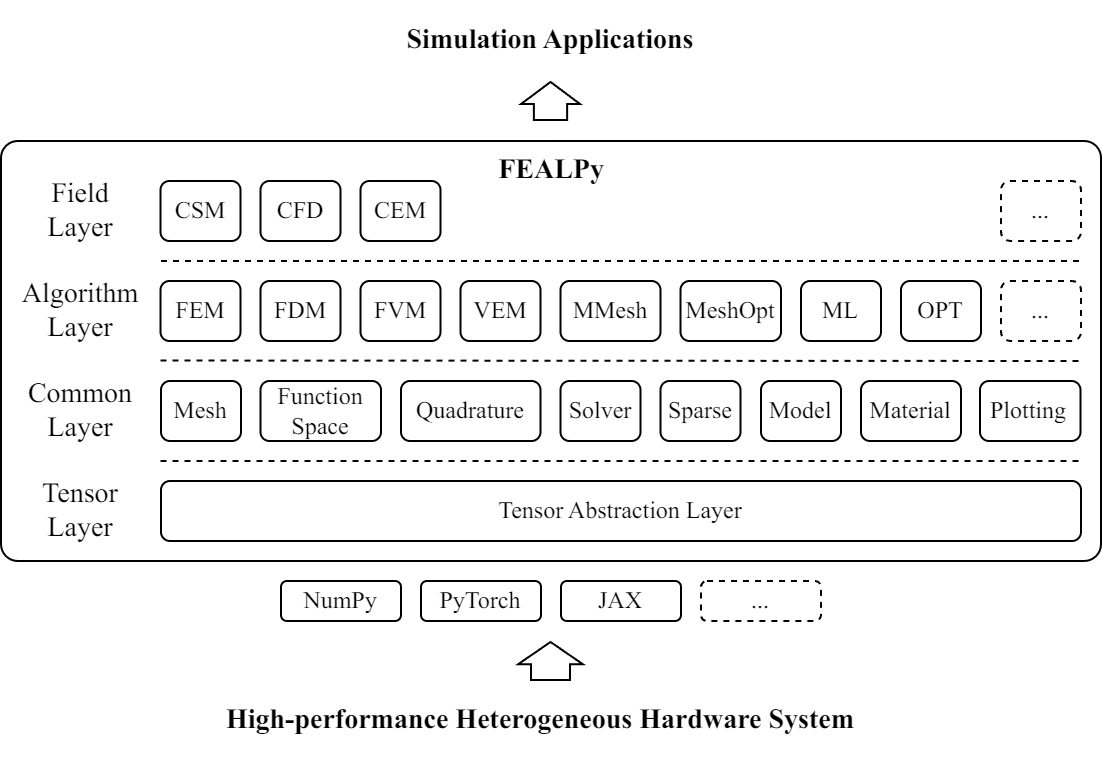}
    \caption{Software architecture of FEALPy}
    \label{fig:fealpy_arch}
\end{figure}

FEALPy is a modular numerical simulation engine built around tensorized
computing.
Its architecture separates mathematical abstractions from computational details
so that numerical algorithms remain independent of specific hardware or
software backends.
Based on the dependencies between tensors, geometry, algorithms and applications,
from bottom to top, we divide FEALPy into four layers: tensor layer,
commonality layer, algorithm layer and field layer, see Figure~\ref{fig:fealpy_arch}.
The tensor layer provides an abstract interface to meet the needs of all upper
modules that are based on tensor computations; the commonality layer provides
basic data structures and algorithms, such as meshes, function spaces, sampling,
integration and solvers; the algorithm layer provides many classic numerical
algorithm frameworks, such as FEM, FDM, FVM, SPH; the field layer enters
each specific field and provides common material settings, physical models,
computational models and examples.

\subsection{Backend Module and Hardware Adaptation}

The backend layer is the foundation of FEALPy's interoperability.
It provides a backend manager (\texttt{BM}) that standardizes tensor operation
interfaces across multiple frameworks.
Every tensor operation in FEALPy—such as array creation, broadcasting,
reduction, or Einstein summation—is routed through this manager, which
dynamically dispatches calls to the active backend.
This design mirrors the multi-backend philosophy adopted by
TensorLy\cite{TENSORLY}, where the same high-level logic can execute on CPUs
or GPUs without code modification.


Backend layer exposes a unified API compatible with the Python Array API
Standard\cite{ARRAYAPI2024}.
It abstracts the complexities of managing different tensor computation backends,
allowing FEALPy to integrate various libraries such as NumPy, PyTorch, JAX,
MindSpore, PaddlePaddle, and others.
This flexibility ensures that users can seamlessly adapt their simulations to
the most suitable tensor backend based on specific hardware requirements.
In addition, as all higher-level modules operate solely on \texttt{BM} to
access APIs, FEALPy ensures that algorithms written once in terms of
abstract tensor operations can be executed on any supported backend.
This not only simplifies development by providing a unified API, but also
ensures long-term maintainability and adaptability as new backends and
hardware architectures emerge.

Key features include:
\begin{itemize}
    \item Unified Data Structure: All computations—from finite element assembly
    to NN operations—are performed on unified tensors, eliminating
    format conversions and device transfers.
    \item Batch Processing: The framework supports mini-batch processing of
    multiple PDE solutions simultaneously, leveraging parallelized computation
    on GPU architectures.
    \item Automatic Differentiation (AD): FEALPy enables gradient computation through
    the entire pipeline, allowing the FEM solver to function as a differentiable
    layer within NNs.
\end{itemize}

This backend abstraction isolates algorithmic logic from low-level performance
tuning, allowing the same solver implementation to benefit from backend-specific
optimizations such as GPU kernels or just-in-time compilation.
It thus forms the computational backbone of FEALPy's multi-paradigm design:
a consistent numerical environment adaptable to both classical PDE solvers and
AI-driven applications.

\subsection{Common Modules with Unified Interfaces}

Common modules such as meshes, function spaces, linear algebra solvers,
samplers, and integrators are implemented based on the tensor abstraction layer.
They defines the essential abstractions shared by different discretization
methods, ensuring that higher-level algorithms operate consistently across
problem types and backends.
The most significant features of these modules are their almost unified
interfaces and strong extensibility, which are particularly emphasized in the
design of FEALPy so that users can have a consistent experience when developing
extensions and higher-level modules.

At its center is the Mesh module, which handle the geometric and
topological structure of computational domains.
The \texttt{Mesh} type is an abstract canonical interface that provides
prescribed geometric and topological information for other modules.
Mesh objects describe nodes, edges, faces, and cells, and store mappings
between these entities.
FEALPy supports both structured and unstructured meshes in 1, 2, and 3
dimensions, along with flexible boundary tagging for applying boundary
conditions.
All mesh data are stored based on tensors, enabling efficient vectorized
operations and direct interaction with the backend layer.




The Function Space module manages function space types, usually based
on meshes, to discrete the solution of problems.
The \texttt{Function Space} type is an abstract interface that defines a set of
functions as the basis, including their values, gradients, and Hessian matrices,
and provides a global index of these functions.
This abstraction allows algorithms such as FEM assembly or gradient evaluation
to be expressed purely through tensor operations, independent of the specific
discretization method.
Different function space subclasses correspond to different function definitions
and global indexing methods, which can serve various types of finite elements,
finite volumes, and even random neural networks, with the same interface,
supporting numerical integration and sampling on functions.


In addition to mesh and function-space management, the common layer includes a set
of modules that serve as reusable building blocks across the framework:
\begin{itemize}
    \item \texttt{Solver}: A variety of algebra solvers for linear systems.
    \item \texttt{Sparse}: A variety of sparse formats (CSR, CSC, and
    COO) supporting multiple backends.
    \item \texttt{Quadrature}: Quadrature rules for numerical
    integration.
    \item \texttt{Plotting}: Lightweight visualization tools for quick
    inspection of meshes and values attached to their entities.
\end{itemize}

Together, these components form a coherent foundation on which discretization
methods can be built.
Their tensor-based design ensures consistency and extensibility:
new numerical schemes can be implemented by reusing these common abstractions
rather than redefining data structures or computation patterns.

\subsection{Basic Algorithms Implemented through Tensor Calculation}

The basic algorithm layer contains a variety of common numerical method modules.
Every computation in this layer, from local element assembly to global matrix
construction, is formulated as tensor operations that can be executed on any
supported backend.
FEM is one of the oldest modules and the starting point for the development of
FEALPy.

The data flow in FEM consists mainly of numerical integration to form element
matrices, global assembly of linear systems, and enforcement of boundary
conditions.
In FEALPy, these processes are organized around three core components:
Integrators, Forms, and Boundary Condition Processors.
\begin{itemize}
    \item Integrators: evaluate local quantities, typically integrals of basis
    functions and their derivatives over mesh entities, and return results as
    element-level tensors. Specialized implementations can also compute fluxes
    across control-volume interfaces, enabling reuse within FVM formulations.
    \item Forms: combine element matrices with global indexing information
    provided by the Function Space module, managing the insertion of nonzero
    entries into sparse matrices.
    \item Boundary Condition Processors: modify assembled systems according to
    Dirichlet constraints.
\end{itemize}
Because all these components operate on tensors, algorithms remain independent
of backends and can automatically benefit from hardware acceleration or
automatic differentiation.
This abstraction also facilitates the extension of existing methods: developers
can prototype new variational forms or discretization strategies simply by
redefining integrators or function spaces without altering the rest of the
framework.

\subsection{The Field Layer: Packaging for Fields}

Built on top of the tensor, common, and algorithmic layers, the field layer
serves as the application interface of FEALPy, where developers and users can
develop specialized subclasses for specific domains.
It organizes domain-specific modules for various areas of computational science
and engineering, such as solid mechanics, fluid dynamics, and electromagnetics.
Each module provides commonly used material models, constitutive relations,
boundary conditions, and example problems tailored to its discipline.

While the present paper focuses on the general software architecture and
tensor-based design of FEALPy, further details of these field-specific
extensions will be discussed in companion works.

\section{Implementation Details and Programming Techniques}
\label{sec:implementation_details}

\subsection{Backend Manager and Adapters}

A significant advantage of FEALPy's Tensor Backend System is its seamless
backend switching capability.
Users can optimize simulations for different computing environments—whether they
are running simulations on CPUs, GPUs, or specialized hardware like TPUs—without
modifying their codebase.
Listing~\ref{lst:switch_backend} snippet demonstrates how FEALPy enables users
to switch between different tensor backends with minimal effort.
This seamless switching capability allows developers to experiment with various
backends and select the one that offers the best performance for their hardware,
ensuring optimized simulation results.
\begin{lstlisting}[language=Python,
    caption={Switching between different tensor backends: By replacing the
    backend name passed to the \\set\_backend function, the object type of the
    output tensor changes. bm here refers to the backend manager, which provides
    standard array APIs.},
    label={lst:switch_backend}
]
from fealpy.backend import bm
from fealpy.mesh import UniformMesh
bm.set_backend("pytorch")
mesh = UniformMesh([0, 1, 0, 1], nx=10, ny=10)
print(type(mesh.entity('cell')))
\end{lstlisting}

The Backend module provides a unified tensor abstraction layer that enables
users to operate different tensor libraries (e.g., NumPy, PyTorch) through a
nearly identical interface.
It consists of two main components: a \texttt{BackendManager} and a set of
backend adapters.
The workflow involves two stages: transmission and translation.

The \texttt{BackendManager} acts as a dispatcher.
When a tensor operation is requested, it routes the call to the currently
selected backend.
If the corresponding adapter has not yet been initialized, it is loaded at that
moment.
This lazy-loading mechanism prevents unnecessary imports and reduces
initialization overhead.

Each adapter implements the Python Array API Standard for a specific
tensor library.
For every interface function required by the standard, the adapter determines
whether and how the target tensor library provides equivalent functionality.
In this context, “satisfied” means that the tensor library already provides an
operation that matches the required functionality of the standard, either
directly or with minor adaptation.
“Not satisfied” means that the tensor library does not provide such
functionality, and the operation must therefore be implemented manually within
the adapter.

Specifically, three situations are considered for each interface function:

\begin{itemize}
    \item \textbf{Directly satisfied}: The tensor library provides a function whose name
    and behavior are consistent with the required interface.
    In this case, the adapter directly references this function.
    \item \textbf{Functionally satisfied but with different naming or signature}:
    The tensor library provides equivalent functionality, but the function name
    or argument format differs from the standard interface.
    In this case, the adapter maps the standard interface function to the
    corresponding library function using a predefined mapping rule.
    \item \textbf{Not satisfied}: The tensor library does not provide a function that
    fulfills the required functionality.
    In this case, the adapter implements the operation explicitly using
    available primitives from the tensor library.
\end{itemize}

Through this mechanism, FEALPy ensures that high-level modules remain
consistent across different backends, while backend-specific differences are
encapsulated within the adapters.

\subsection{Mesh Data Structure Based on Tensors}

Based on the tensor abstraction interface, we designed the mesh data structure
in FEALPy carefully.
The Cartesian coordinate of points are recorded in a floating-point tensor,
which is then combined into higher-dimensional shapes through a series of
non-negative integer tensors.
We adopt the notation shown in Table~\ref{tab:notations}.
\begin{table}[htbp]
\centering
\begin{tabular}{cl}
\hline
\textbf{Notation} & \multicolumn{1}{c}{\textbf{Description}} \\ \hline
GD                & Geometry dimensions                      \\
TD                & Topological dimensions                   \\
NN                & Number of nodes                          \\
NE                & Number of edges                          \\
NF                & Number of faces                          \\
NC                & Number of cells                          \\ \hline
\end{tabular}
\caption{Common notations in FEALPy. Node, edge, face and cell are mesh
entities.}
\label{tab:notations}
\end{table}

\paragraph{Nodes.}
The position of a node is determined by a vector of length \texttt{GD}, so all
nodes in a mesh is stored as a tensor with shape \texttt{[NN, GD]}.
Once the tensor is stored, not only the node positions are identified,
but also their global indices.

\paragraph{Line segments, triangles and quadrangles.}
Two vertices can form a line segment by connecting them in a directed way, so
a 2-tuple containing the global indices of the start and end points can
represent an edge.
Edges of a mesh can be represented as a tensor with shape \texttt{[NE, 2]} in
almost all unstructured mesh types in FEALPy.
2D triangles are constructed in FEALPy by connecting three vertices in
counter-clockwise order, thus \texttt{[num\_triangles, 3]} records a
collection of them in a mesh.
Listing~\ref{lst:mesh_init} is an example for a mesh containing 2 triangles
divided from a square, and the rules for other shapes are similar.
\begin{lstlisting}[language=Python,
    caption={A triangular mesh is stored as two tensors—\texttt{node} and \texttt{cell}},
    label={lst:mesh_init}
]
from fealpy.mesh import TriangleMesh
node = bm.tensor([
    [0. 0.] # 0
    [0. 1.] # 1
    [1. 0.] # 2
    [1. 1.] # 3
], dtype=bm.float64)
cell = bm.tensor([
    [2 3 0] # 0
    [1 0 3] # 1
], dtype=bm.int64)
mesh = TriangleMesh(node, cell)
\end{lstlisting}

\begin{figure*}[htbp]
    \centering
    \begin{minipage}{0.50\textwidth}
        \includegraphics[width=1.0\textwidth]{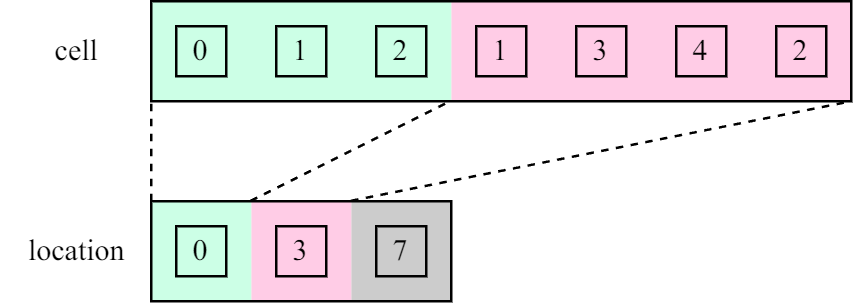}
        \caption{For polygons, an additional pointer tensor is used to record
        the starting index of each polygon.}
        \label{fig:polygon}
    \end{minipage}
    \begin{minipage}{0.30\textwidth}
        \includegraphics[width=1.0\textwidth]{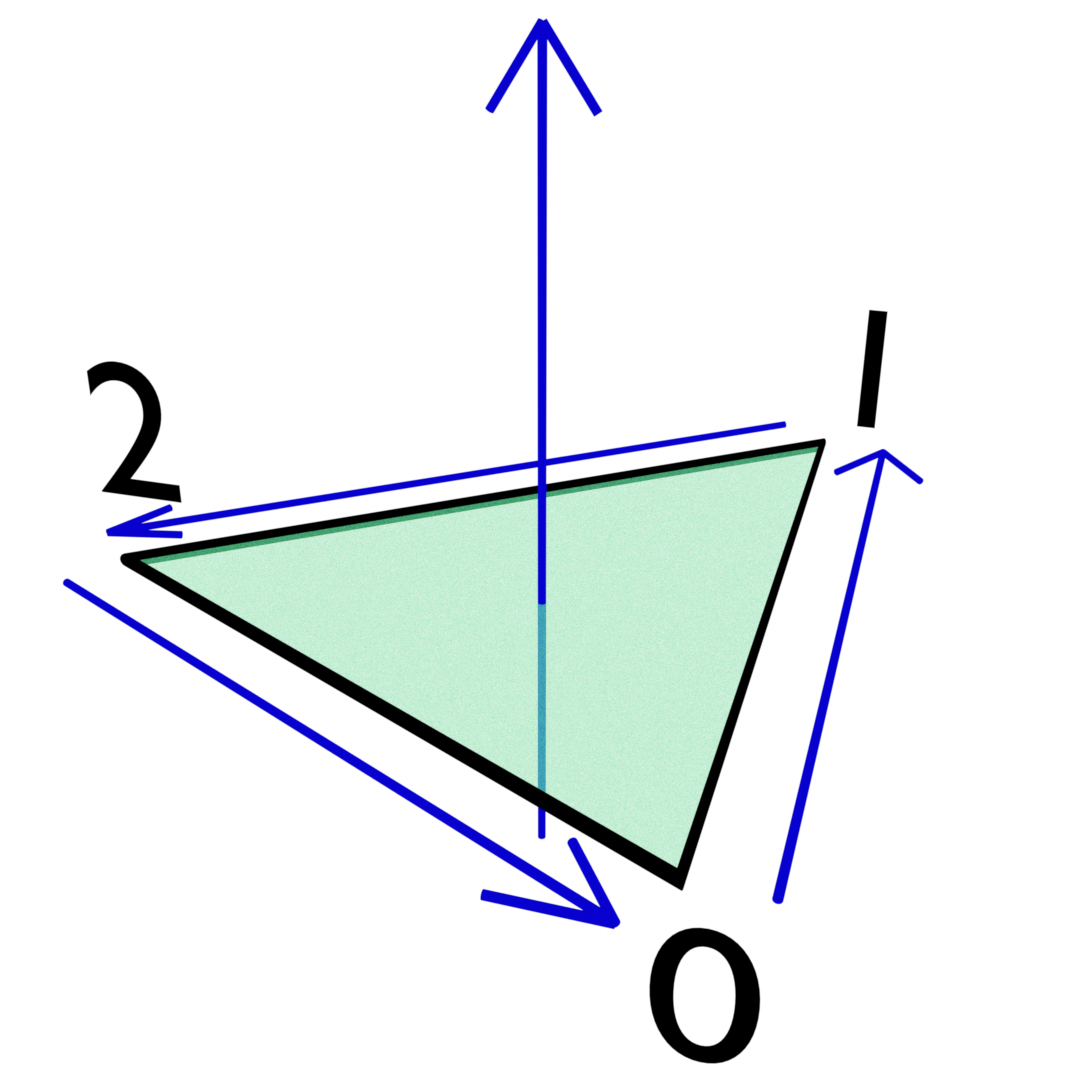}
        \caption{Orientation of a 3D triangle is defined by the right-hand rule.}
        \label{fig:orientation}
    \end{minipage}
\end{figure*}

\paragraph{Polygons.}
Polygons describe a group of shapes with inconsistent numbers of edges, and
these shapes cannot be combined to form 2D tensors similar to the rules for
triangles and quadrangles.
All vertices from multiple polygons are directly arranged in a 1D tensor, and a
pointer tensor is used to record the starting position (with an additional total
number of vertices as the last element) of each polygon, as shown in Figure~
\ref{fig:polygon}.
The length of this pointer tensor is exactly \texttt{num\_polygons + 1}, like the
\texttt{indptr} in CSR format sparse matrices.

\paragraph{Orientation.}
In the case of unstructured meshes, FEALPy always forms these topological
2-dimensional shapes by ordering the vertices in a circular manner, so it is
natural to define the orientation by the ordering direction.
The right-hand rule is used to define the orientation in FEALPy.
See Figure~\ref{fig:orientation}.


\paragraph{Mesh entities.}
FEALPy is more inclined to regard these organizations of nodes as entities by
their topological dimensions, instead of their shapes.
Nodes—which only have positional information—are 0-dimensional entities,
while edges are 1-dimensional entities.
Cells refer to entities with the highest topological dimension,
such as any triangle in a triangular mesh.
Faces refer to entities which is 1 less dimensional than cells, and they usually
locate on mesh boundaries or interfaces between cells.
The topological dimension of cells (i.e. the highest topological dimension) is
also the topological dimension of the mesh.
In particular, edges and faces refer to the same kind of entities in a
topologically 2D mesh.

\paragraph{Topological relations.}
Some topological relationships between mesh entities are available in FEALPy,
such as \texttt{face\_to\_cell} in the shape of \texttt{[NF, 4]}.
A deduplication process is needed for constructing these relationships, because
lower-dimensional entities, e.g. edges and faces, are bound to have duplicates
if they are simply collected from all higher-dimensional cells.
Table~\ref{tab:face_to_cell} shows what the four columns of
\texttt{face\_to\_cell} mean.
FEALPy deduplicates the lower-dimensional entities using \texttt{lexsort},
through which the first and the last occurrences of the same faces are
extracted from the total.
Since the right-hand rule is used for 2D face orientation, the cell where
a face firstly appears is the cell the face faces to;
and since the counter-clockwise rule is used for 1D edge orientation, the cell
where an edge firstly appears is the cell on the left side of the edge.
\begin{table}[htbp]
\centering
\begin{tabular}{cccc}
\hline
\textbf{0}                                                         & \textbf{1}                                                         & \textbf{2}                                                         & \textbf{3}                                                         \\ \hline
\begin{tabular}[c]{@{}c@{}}cells indices\\ left/front\end{tabular} & \begin{tabular}[c]{@{}c@{}}cells indices\\ right/back\end{tabular} & \begin{tabular}[c]{@{}c@{}}local indices\\ left/front\end{tabular} & \begin{tabular}[c]{@{}c@{}}local indices\\ right/back\end{tabular} \\ \hline
\end{tabular}
\caption{Definition of the 4 columns in \texttt{face\_to\_cell}.
Each row corresponds to a face after deduplication, and the first two columns
represent the indices of the two cells that this face belongs to, while the last
two columns represent the local indices of this face in these two cells.}
\label{tab:face_to_cell}
\end{table}

\paragraph{Mesh classes.}
By grouping various entities together and adding related methods, a mesh class
can be constructed.
All supported mesh types are listed in Table~\ref{tab:mesh_types}.
The Mesh module is designed to maximize the reuse of methods through an
reasonable inheritance structure and keep the interface of methods consistent
between different mesh types.
Therefore, meshes can be simply replaced without the need to adapt the rest of
the script.

Listing~\ref{lst:mesh_class} shows mesh types that currently supported in
FEALPy.
And we would like to make some extra explanations here:
\begin{itemize}
    \item The \texttt{NodeMesh} is a special mesh type that only contains nodes
    without any higher-dimensional entities, and the topological relationships
    between nodes are dynamically constructed, for example, by connecting
    adjacent nodes leveraging a KD-tree.
    This is particularly useful for methods like smooth particle hydrodynamics
    (SPH) and deep learning methods that operate on point clouds.
    \item Each edge of a mesh can be divided into a pair of opposite edges if
    the directionality is considered, where each in the pair is called a
    'half edge'.
    Therefore, the entire mesh, as a \texttt{HalfEdgeMesh}, has a total of
    $2\texttt{NE}$ half edges.

\end{itemize}
\begin{table}[htbp]
\centering
\begin{tabular}{ll|ll|ll}
\hline
\textbf{Cell Type}   & \textbf{TD} & \textbf{Cell Type}  & \textbf{TD} & \textbf{Mesh Type} & \textbf{TD} \\ \hline
Interval             & 1           & Edge                & 1           & Node Mesh          & 0           \\
Quadrangle           & 2           & Triangle            & 2           & HalfEdge Mesh      & 1           \\
Tetrahedron          & 3           & Polygon             & 2           & Dart Mesh          & 2/3         \\
Prism                & 3           & Hexahedron          & 3           & Uniform Mesh       & 1/2/3       \\
Lagrange Triangle    & 2           & Lagrange Quadrangle & 2           &                    &             \\
Lagrange Tetrahedron & 3           & Lagrange Hexahedron & 3           &                    &             \\ \hline
\end{tabular}
\caption{Types of unstructured meshes (left and middle columns) and others
         (right column) in FEALPy and their supported TDs.}
\label{tab:mesh_types}
\end{table}

\begin{lstlisting}[language=Python,
    caption={The mesh type can be simply replaced directly. The code shows
    the replacement from TriangleMesh to QuadrangleMesh.},
    label={lst:mesh_class}
]
from fealpy.mesh import TriangleMesh, QuadrangleMesh
from fealpy.functionspace import LagrangeFESpace

if meshtype == "TRIA":
    mesh = TriangleMesh.from_box()
elif meshtype == "QUAD":
    mesh = QuadrangleMesh.from_box()
else:
    raise ValueError()

space = LagrangeFESpace(mesh, p=1)
...
\end{lstlisting}

\subsection{Basis and DoF Management}

Numbering the global degrees-of-freedom (DoF) and inferring their relationships
with local DoFs in mesh entities is called DoF management in FEALPy.

\begin{figure}[htbp]
    \centering
    \includegraphics[width=0.50\textwidth]{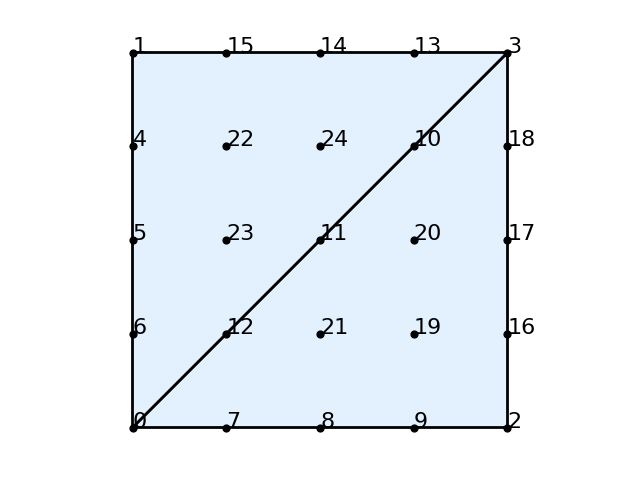}
    \caption{Example of global DoF numbering ($p=4$) in a triangle in FEALPy.}
    \label{fig:dof_numbering}
\end{figure}

For most finite element methods or their variants, DoFs are arranged by
interpolation points on the mesh entities, and there are two main types of DoF
management provided in FEALPy in this way.
The first type is continuous, where adjacent mesh cells may share a part of DoFs
at their interface.
Section 6 of \cite{CHEN2023GEO} has discussed the global indexing rules
for Lagrange interpolation points in detail, which is actually the main idea
behind FEALPy's implementation of various finite element DoF management.
In order to provide a non-repeating numbering that is as consistent as possible
for all mesh entities, we adopt the order of DoFs from low-dimensional entities
to high-dimensional entities in the order of "node-edge-face-cell", as shown in
Figure~\ref{fig:dof_numbering}.
When numbering the DoFs inside two adjacent cells, the DoFs at their interface
are already numbered.
The second type is discontinuous, where adjacent mesh cells do not share any
DoF.
Non-repeating numbers can be generated directly from the correspondence between
cells and their local DoFs.

Take the \texttt{Lagrange} space as an example—both basis generation and DoF
management in 1, 2, and 3 dimensions can be supported in FEALPy with arbitrary
order.
Because they are both constructed from the multiple indices on mesh entities,
which are generated as a multi-dimensional tensor, independent
of the spatial dimensions and the order of the indices.


Based on different basis function implementations and mappings of local-global
DoFs, FEALPy already supports a variety of function spaces summarized in Table~
\ref{tab:function_spaces}.
\begin{table}[htbp]
\centering
\begin{tabular}{llc}
\hline
\textbf{Element Types}  & \textbf{GD} & \textbf{References}                      \\ \hline
Lagrange                & Any         &                                          \\
Bernstein               & Any         & \cite{KIRBY2011-Bernstein}               \\
$C^m$ Conforming        & 2, 3        & \cite{CHENLONG2021, HUJUN2023, CHEN2025} \\
Hu-Zhang                & 2, 3        & \cite{HU2015-HuZhang}                    \\
First Nédélec           & 2, 3        & \cite{NEDELEC1986}                       \\
Second Nédélec          & 2, 3        & \cite{NEDELEC1986}, \cite{CHEN2023GEO}   \\
Raviart-Thomas          & 2, 3        & \cite{RAVIART1977}                       \\
Brezzi-Douglas-Marini   & 2, 3        & \cite{BREZZI1985-BDM}, \cite{CHEN2023GEO}\\
Parametric Lagrange     & Any         &                                          \\ \hline
\end{tabular}
\caption{Element types in FEALPy and their supported GDs.}
\label{tab:function_spaces}
\end{table}

\paragraph{Tensor function spaces.}
Given a scalar function space, the \texttt{TensorFunctionSpace} class in
FEALPy constructs a tensor-valued function space by assigning a prescribed
tensor shape to each scalar degree of freedom (DoF).
A leading or trailing `-1` in the shape specifies the ordering of the tensor DoFs:
\begin{itemize}
    \item A leading `-1` indicates that the scalar DoF index serves as the outer index,
    so that the tensor component index varies fastest.
    \item A trailing `-1` indicates that the scalar DoF index serves as the inner index,
    so that the scalar DoF index varies fastest within each tensor component.
\end{itemize}
Listing~\ref{lst:tensor_function_space} illustrates the construction of a
tensor function space with shape \texttt{(-1, 2)}—a 2-component vector-valued
function space where the scalar DoF index is the outer index.

\begin{lstlisting}[language=Python,
    caption={A tensor-valued function space was constructed by shape (-1, 2),
    from a given scalar space.
    The leading -1 means that the scalar DoF index serves as the outer index,
    so that the tensor component index varies fastest.},
    label={lst:tensor_function_space}
]
from fealpy.functionspace import (
    LagrangeFESpace, TensorFunctionSpace
)
...
scalar_space = LagrangeFESpace(mesh=mesh, p=1)
tensor_space = TensorFunctionSpace(scalar_space, shape=(-1, 2))
\end{lstlisting}

\subsection{Sampling and Quadrature}

Sampling and quadrature are also essential tools for many numerical methods.
The \texttt{Quadrature} module stores precomputed barycentric coordinates of
quadrature points and their weights generated from symmetric quadrature rules
on simplexes \cite{WILLIAMS201418}, including line segments, triangles, and
tetrahedrons.
While quadrature points on $d$-hypercubes, such as quadrangles and hexahedrons,
are generated via tensor products of 1D Gauss-Legendre quadrature rules
(Listing~\ref{lst:quadrature_tensor_product}).
Similarly, triangular prisms and pyramids can also get their quadrature rules
via tensor products.
\begin{lstlisting}[language=Python,
    caption={Quadrature points on quadrangles can be generated via tensor
    products of two 1D Gauss-Legendre quadrature rules.},
    label={lst:quadrature_tensor_product}
]
quad_pts = (quad_pts_1d_a, quad_pts_1d_b)
weights = bm.einsum(
    'i,j->ij', weights_1d_a, weights_1d_b
).reshape(-1)
\end{lstlisting}

Samplers are designed for generating equidistant or random collocation points
for deep learning methods.
For instance, \texttt{ISampler} and \texttt{BoxBoundarySampler} can generate
sampling points in the interior and on the boundary of a given rectangular
region respectively;
while \texttt{MeshSampler} can draw sampling points on given homogeneous mesh
entities.


\subsection{Sparse Format and Matrix Assembly}

FEALPy implements sparse matrix representations in COO and CSR formats.
These implementations are backend-agnostic: all sparse data structures are
constructed entirely from the tensor abstraction layer provided by the
Backend module, and therefore operate uniformly across all supported tensor
libraries.

\paragraph{Extended sparse matrix formats.}
Let $d$ denote the number of tensor axes and $\mathrm{nnz}$ the number of
nonzero entries.
In addition to standard COO and CSR definitions, FEALPy extends these formats
in two aspects.

First, the \texttt{data} field may be \texttt{None}, which corresponds to
structure-only sparse matrices (e.g., Boolean adjacency matrices), where
only index information is required.

\begin{figure}
    \centering
    \includegraphics[width=0.4\textwidth]{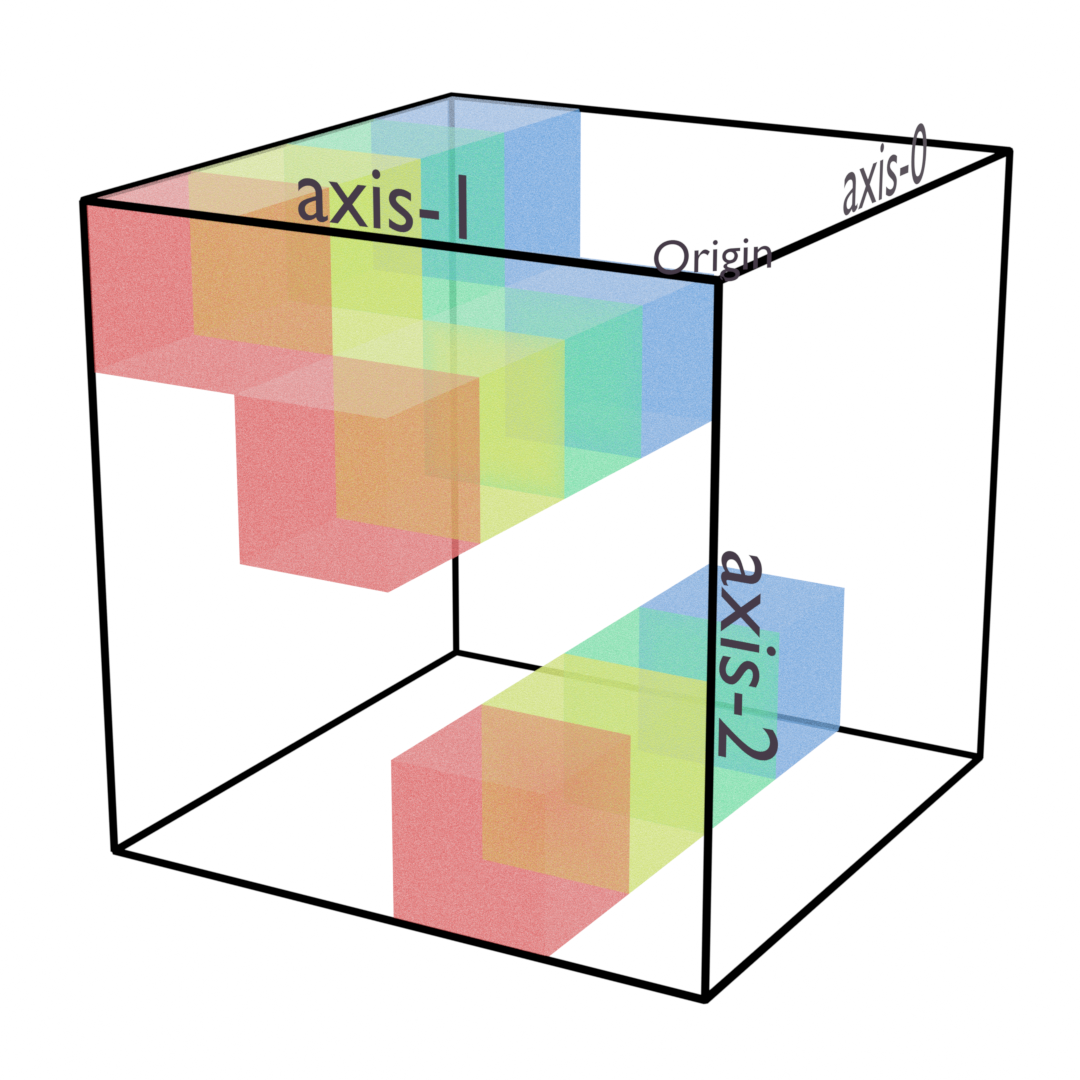}
    \caption{Hybrid COO: a batched sparse matrix containing many sparse matrices
    with the same non-zero pattern. This is useful in assembling matrices for
    batches of equations with different coefficients in FEM.}
    \label{fig:HCOO}
\end{figure}

Second, the \texttt{data} tensor may contain a leading batch axis.
Specifically, if the sparse tensor represents a batch of matrices sharing
the same sparsity pattern, the index tensor has shape $(d, \mathrm{nnz})$,
while the data tensor has shape $(b, \mathrm{nnz})$, where $b$ is the batch
size.
This hybrid representation (Figure~\ref{fig:HCOO}) enables efficient
storage of multiple sparse matrices with identical nonzero structure but
different numerical values.
Such a design is particularly useful in finite element computations where
multiple linear systems with varying coefficients share the same connectivity
pattern.

\paragraph{Batched matrix assembly.}
In the finite element module, global sparse matrices are assembled by
combining local element matrices produced by the \texttt{Integrator} with
degree-of-freedom (DoF) mappings provided by the \texttt{FunctionSpace}.
The \texttt{Forms} component coordinates this process using the Sparse module.

Let $\texttt{cell\_matrix}$ denote the local stiffness matrices with shape
$$
(\texttt{num\_cells}, ~\texttt{num\_quadpts}, ~\texttt{num\_basis}, ~\texttt{num\_basis}),
$$
and let $\texttt{ue2dof}$ and $\texttt{ve2dof}$ be the DoF mappings from
local basis functions to global indices, both of shape
$(\texttt{num\_cells}, ~\texttt{num\_dofs}).$
The global row and column indices of all nonzero entries are generated
by broadcasting the DoF mappings to match the local tensor shape.
Listing~\ref{lst:global_assembly} illustrates this procedure.

\begin{lstlisting}[language=Python,
    caption={Assembly of global stiffness matrices. The local\_shape equals to
    (num\_cells, num\_quadpts, num\_basis, num\_basis) which stores the shape
    of cell\_matrix; ue2dof and ve2dof are DoF maps, in the shape of
    (num\_cells, num\_dofs).},
    label={lst:global_assembly}
]
I = bm.broadcast_to(ve2dof[:, :, None], local_shape)
J = bm.broadcast_to(ue2dof[:, None, :], local_shape)
indices = bm.stack([I.ravel(), J.ravel()], axis=0)
cell_matrix = bm.reshape(cell_matrix, (-1,))
# The shape will be (batch_size, -1) for batched assembly.
M = M.add(COOTensor(indices, cell_matrix, sparse_shape))
\end{lstlisting}

The broadcasted tensors \texttt{I} and \texttt{J} represent the global row
and column indices of each local contribution.
After flattening, these indices are stacked to form the COO \texttt{index}
tensor of shape $(2, \mathrm{nnz})$, where $\mathrm{nnz}$ equals the total
number of local contributions.
The local matrices are reshaped accordingly to produce the associated COO
\texttt{data} tensor.
The global sparse matrix is then constructed by accumulating these contributions
via the \texttt{add} operation.

For batched problems, the local contribution tensor acquires an additional leading
dimension corresponding to the batch size, resulting in a sparse tensor
with shape $(b, \mathrm{nnz})$ in the \texttt{data} field.

\paragraph{Current limitations.}
At present, mainstream tensor computation libraries provide limited
low-level (C/C++) support for batched sparse matrix–vector multiplication
when matrices share identical sparsity patterns.
Although correct implementations can be constructed at the Python level using
standard tensor operations, these approaches typically rely on repeated sparse
operations or indirect indexing, which introduces additional overhead.
As a result, solver performance may degrade in large-scale batched
computations.

Improving backend-level support for batched sparse linear algebra remains
an important direction for future development.
And note that our unified tensor abstraction layer allows us to integrate future
improvements or user C/C++ extensions
(by updating the existing adapters or adding a new backend adapter)
without changing the high-level algorithm code.
So the current framework will not prevent future improvements in this area.

\section{Numerical Examples}

\subsection{3D Linear Elasticity Problem with Lagrange Finite Elements}
This subsection presents a fundamental numerical example in solid mechanics,
designed to verify the basic correctness and computational accuracy of FEALPy.
Given computational domain $\Omega=[0,1]^3$, the material is modeled as a
homogeneous, isotropic, linear elastic solid with Lamé constants $\lambda=1$
and shear modulus $\mu=1$.
The governing system of equations is given by the equilibrium equation,
kinematic strain-displacement relation, and the constitutive law
(Generalized Hooke's Law):
\begin{align}
-\nabla \cdot \bsigma &= \bfb \quad \text{in } \Omega \\
\bvarepsilon &= \frac{1}{2} (\nabla \bfu + (\nabla \bfu)^T) \\
\bsigma &= \lambda \, \mathrm{tr}(\bvarepsilon) \bfI + 2\mu \bvarepsilon
\end{align}
where $\bfu$ is the displacement vector, $\bvarepsilon$ is the strain tensor,
and $\bsigma$ is the Cauchy stress tensor.
The exact displacement field $\bfu = (u_1, u_2, u_3)$ is chosen to be:
$$
\bfu = \begin{bmatrix}
200\mu(x-x^2)^2(2y^3-3y^2+y)(2z^3-3z^2+z)\\
-100\mu(y-y^2)^2(2x^3-3x^2+x)(2z^3-3z^2+z)\\
-100\mu(z-z^2)^2(2y^3-3y^2+y)(2z^3-3z^2+z)
\end{bmatrix}
$$
The corresponding body force density $\bfb$, which satisfies the equilibrium
equation above, is derived accordingly and is:
$$
\bfb = \begin{bmatrix}
-400\mu(2y-1)(2z-1)[3(x^2-x)^2(y^2-y+z^2-z)+(1-6x+6x^2)(y^2-y)(z^2-z)]\\
200\mu(2x-1)(2z-1)[3(y^2-y)^2(x^2-x+z^2-z)+(1-6y+6y^2)(x^2-x)(z^2-z)]\\
200\mu(2x-1)(2y-1)[3(z^2-z)^2(x^2-x+y^2-y)+(1-6z+6z^2)(x^2-x)(y^2-y)]
\end{bmatrix}
$$
A non-homogeneous Dirichlet boundary condition is applied over the entire
boundary $\partial \Omega$.

In discretization, the domain $\Omega$ was initially partitioned into a
uniform $4 \times 4 \times 4$ mesh of tetrahedrons, and was successively
uniformly refined four times.
Finite element spaces $V_h \subset H^1(\Omega)$ was constructed using linear or
quadratic Lagrange elements on each mesh level.
Listing~\ref{lst:linear_elasticity} shows the construction of the finite
element space.
\begin{lstlisting}[language=Python,
    caption={A tensor function space was constructed on a tetrahedral mesh.
    The 'shape' parameter specifies the tensor shape and DoF arrangement.},
    label={lst:linear_elasticity}
]
from fealpy.mesh import TetrahedronMesh
from fealpy.functionspace import (
    LagrangeFESpace, TensorFunctionSpace
)
mesh = TetrahedronMesh.from_box(box=domain, nx=4, ny=4, nz=4)
space = LagrangeFESpace(mesh=mesh, p=1, ctype='C')
tensor_space = TensorFunctionSpace(space, shape=(-1, 2))
\end{lstlisting}

The weak formulation of the problem is: Find $\bfu_h \in V_h$ such that
$$
\int_\Omega \bsigma(\bfu_h) : \bvarepsilon(\bfv_h) \rmd\Omega
= \int_\Omega \bfb \cdot \bfv_h \, d\Omega \quad \forall \bfv_h \in V_h^0
$$
where $V_h^0$ is the test space satisfying homogeneous Dirichlet boundary
conditions.
See Listing~\ref{lst:linear_elasticity_form} for the code.
The non-homogeneous Dirichlet conditions are enforced in a strong form.
\begin{lstlisting}[language=Python,
    caption={Construct a bilinear form based on the space, and integrators
    can be added to the form to define the weak form. The assembly method
    returns a sparse matrix. The construction of linear forms is similar.},
    label={lst:linear_elasticity_form}
]
from fealpy.fem import BilinearForm, LinearElasticIntegrator
from fealpy.material import IsotropicLinearElasticMaterial

material = IsotropicLinearElasticMaterial(
    lame_lambda=1.0,
    shear_modulus=1.0,
    plane_type='plane_strain',
)
bform = BilinearForm(tensor_space)
bform.add_integrator(
    LinearElasticIntegrator(material=material)
)
A = bform.assembly()
\end{lstlisting}

The accuracy of FEALPy's computation was quantified by calculating the
$L^2$-norm error of the displacement:
$$\text{Error} = \Vert \bfu_{\text{exact}} - \bfu_h \Vert_{L^2(\Omega)}$$
The computed $L^2$ errors, order of convergence and the corresponding mesh sizes
are summarized in Table~\ref{tab:LinearElasticity} and
Table~\ref{tab:LinearElasticity_p10}, where the linear systems were solved by Pardiso.
\begin{table}[htbp]
\centering
\begin{tabular}{cccccccc}
\hline
\textbf{p} & \textbf{h}  & \textbf{1/2} & \textbf{1/4} & \textbf{1/8}    & \textbf{1/16}   & \textbf{1/32}   & \textbf{1/64}   \\ \hline
1          & DoF         & -            & 375          & 2187            & 14739           & 107811          & 823875          \\
           & $L^2$ Error & -            & 3.0271e-2    & 1.1616e-2       & 3.5874e-3       & 9.8912e-4       & 2.5827e-4       \\
           & Order       & -            & -            & 1.3818          & 1.6951          & 1.8587          & \textbf{1.9372} \\ \hline
2          & DoF         & -            & 2187         & 14739           & 107811          & 823875          &                 \\
           & $L^2$ Error & -            & 5.0775e-3    & 6.0218e-4       & 6.8054e-5       & 7.9984e-6       &                 \\
           & Order       & -            & -            & 3.0759          & 3.1454          & \textbf{3.0889} &                 \\ \hline
4          & DoF         & 2187         & 14739        & 107811          & 823875          &                 &                 \\
           & $L^2$ Error & 1.3435e-03   & 6.0182e-05   & 2.1914e-06      & 7.2287e-08      &                 &                 \\
           & Order       & -            & 4.4805       & 4.7794          & \textbf{4.9219} &                 &                 \\ \hline
8          & DoF         & 14739        & 107811       & 823875          &                 &                 &                 \\
           & $L^2$ Error & 2.2096e-06   & 4.9601e-09   & 9.9604e-12      &                 &                 &                 \\
           & Order       & -            & 8.7992       & \textbf{8.9600} &                 &                 &                 \\ \hline
\end{tabular}
\caption{Convergence results for the linear elasticity problem.
The $L^2$-norm error decreases systematically as the mesh is refined.}
\label{tab:LinearElasticity}
\end{table}

\begin{table}[htbp]
\centering
\begin{tabular}{ccccc}
\hline
\textbf{p} & \textbf{h}  & \textbf{1} & \textbf{1/2} & \textbf{1/4} \\ \hline
10         & DoF         & 3993       & 27783        & 206763       \\
           & $L^2$ Error & 1.5787e-14 & 6.5569e-14   & 9.3363e-14   \\
           & Order       & -          & -2.0543      & -0.5098      \\ \hline
\end{tabular}
\caption{Convergence results for the linear elasticity problem with $p=10$.
This is an exact solution up to machine precision and the error of
linear system solver.}
\label{tab:LinearElasticity_p10}
\end{table}

The calculated convergence order approaches the optimal theoretical order,
(e.g. 2 for linear elements, and 3 for quadratic elements in the $L^2$-norm)
strongly validating the correct implementation of the FEM for
linear elasticity within FEALPy and confirms its numerical precision.

\subsection{3D Biharmonic Equation with Smooth Finite Elements}
The biharmonic equation arises in models of thin plate deformation and Stokes
flow problems. Here, the following problem is considered on a unit cube domain
$\Omega=(0, 1)^3$:
$$\Delta^2u=f$$
with manufactured solution $u_{\text{exact}}=\sin(4x)\sin(4y)\sin(4z)$,
and Dirichlet boundary conditions prescribing
both $u$ and $\partial_nu$ on $\partial\Omega$.
The weak formulation requires $H^2$-conforming (i.e., $C^1$-continuous) finite
element spaces, which are non-trivial to construct on simplicial meshes in 3D.

We use the $C^1$-conforming finite element space proposed by
\cite{CHENLONG2021, HUJUN2023}, with shape functions of degree $k\ge 9$
on tetrahedral meshes.

Benefiting from FEALPy's modular design, users are allowed to
implement new finite elements by defining only the local basis functions and
degrees of freedom, while reusing existing infrastructure for mesh handling,
assembly, and linear algebra.
To solve this problem, we simply replaced the standard Lagrange
element with the smooth element ($C^m$ Conforming space, implemented by
\cite{CHEN2025}), on a initially
$1 \times 1 \times 1$ uniform mesh of tetrahedrons.
See Listing~\ref{lst:smooth_finite_element}.
\begin{lstlisting}[language=Python,
    caption={The Smooth Finite Element can be implemented by introducing a new
    function space and integrator.},
    label={lst:smooth_finite_element}
]
...
mesh = TriangleMesh.from_box([0, 1, 0, 1], n, n)
space = CmConformingFESpace2d(mesh, p, 1, isCornerNode)
bform = BilinearForm(space)
MLI = MthLaplaceIntegrator(m=2, coef=1, q=p+4)
bform.add_integrator(MLI)
A = bform.assembly()
...
\end{lstlisting}

Table~\ref{tab:Biharmonic} shows the convergence on a sequence of uniformly
refined tetrahedral meshes, confirming optimal convergence rates in the
$L^2$-norm, as theoretically expected.
\begin{table}[htbp]
\centering
\begin{tabular}{ccccc}
\hline
\textbf{h}                                                          & 1         & 1/2       & 1/4        & 1/8        \\ \hline
\textbf{DoF}                                                        & 582       & 2761      & 16791      & 116971     \\ \hline
$\Vert u_{\text{exact}}-u_h\Vert_{L^2}$                             & 6.6494e-2 & 6.0048e-5 & 7.1664e-8  & 5.4160e-11 \\
Order                                                               & -         & 10.1129   & 9.7106     & 10.3698    \\ \hline
$\Vert\nabla u_{\text{exact}} - \nabla u_h\Vert_{L^2(\Omega)}$      & 4.1421e-1 & 1.1596e-3 & 2.8335e-6  & 4.5102e-9  \\
Order                                                               & -         & 8.4806    & 8.6768     & 9.2952     \\ \hline
$\Vert\nabla^2 u_{\text{exact}} - \nabla^2 u_h\Vert_{L^2(\Omega)}$  & 3.4483e+0 & 2.1847e-2 & 1.0433e-4  & 3.4023e-7  \\
Order                                                               & -         & 7.3023    & 7.7101     & 8.2605     \\ \hline
\end{tabular}
\caption{Convergence results for the 3D biharmonic equation with $C^1$-tet elements,
and shape functions of degree $k=9$.}
\label{tab:Biharmonic}
\end{table}

This example demonstrates the flexibility of FEALPy in implementing non-standard
finite elements.
By leveraging a modular architecture, new elements
(such as the $C^1$-conforming tetrahedral element) can be integrated
seamlessly — users need only provide the local basis and degree-of-freedom
definitions, while all other components
(e.g., assembly, boundary conditions, solvers) remain reusable.

\subsection{Burgers' Equation with Moving Mesh}

The moving mesh method, also known as the r-adaptive method, relocates a fixed
number of mesh nodes according to user-defined monitor functions.
It has been widely applied in fields such as fluid dynamics, heat transfer,
and combustion simulation, where solutions often exhibit sharp gradients or
localized features, such as shock waves or boundary layers, amid relatively
smooth regions.
While fixed uniform meshes may fail to capture such features efficiently without
excessive refinement, the moving mesh method adaptively concentrates mesh nodes
in regions of high variation, significantly improving computational accuracy
and efficiency without altering the mesh topology.

In FEALPy, we have implemented a high-dimensional, vectorized moving mesh
algorithm based on the Harmonic map\cite{RUOLI2001, RUOLI2002} and the MMPDE
gradient flow approach\cite{HUANG2001}.
A \texttt{MMesher} class is provided to control the moving mesh process, as
presented in Listing~\ref{lst:moving_mesh}.
\begin{lstlisting}[language=Python,
    caption={MMesher is designed to control the moving mesh process. It is
    irrelevant to mesh types and FE spaces.},
    label={lst:moving_mesh}
]
from fealpy.mmesh.mmesher import MMesher
...
mm = MMesher(mesh=mesh, uh=uh, space=space, beta=beta,
             is_multi_phy=False, config=config)
mm.initialize()
mm.set_interpolation_method('linear')
mm.set_monitor('arc_length')
mm.set_mol_method('projector')
if t == 0.0: # preprocess
    mm.set_interpolation_method('solution')
    mm.run()
    mm.set_interpolation_method('linear')
mesh, uh = mm.run()
uh[:] = ... # finish assembly and solve
mm.instance.uh = uh # update the instance
\end{lstlisting}
The implementation supports both NumPy and PyTorch backends and can run on CPUs
and GPUs.
It is dimension-independent and allows seamless switching between mesh types,
including high-order curved meshes on planes.

Let's consider a scalar Burgers' equation on $\Omega = [0,1]^2$ with time
$T = [0,2]$
\begin{equation}
\left\{\begin{aligned}
\frac{\partial u}{\partial t} + u(u_x+ u_y) &= a\Delta u, \quad&\bfx \in \Omega, t\in T, \\
u(\bfx , 0) &= u_0(\bfx), \quad&\bfx \in \Omega, \\
u(\bfx , t) &= g(\bfx , t), \quad&\bfx \in \partial \Omega, t \in T.
\end{aligned}\right.
\end{equation}
to demonstrate the correctness, efficiency, and precision of the implementation.
The exact solution is given by
$$
u(x,y,t; a) = \frac{1}{1+ \exp((x+y - t)/2a)}, \quad \text{on }\Omega\times T, a > 0.
$$
where $a>0$ is the viscosity coefficient.
This solution features a sharp transition layer along the line $x+y=t$,
propagating diagonally across the domain.

\begin{figure*}[htbp]
    \centering
    \begin{minipage}{0.4\textwidth}
        \includegraphics[width=1.0\textwidth]{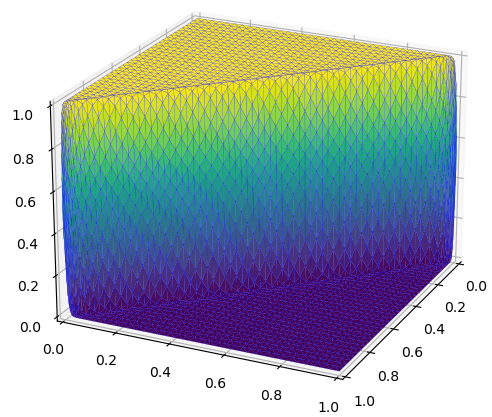}
    \end{minipage}
    \begin{minipage}{0.4\textwidth}
        \includegraphics[width=1.0\textwidth]{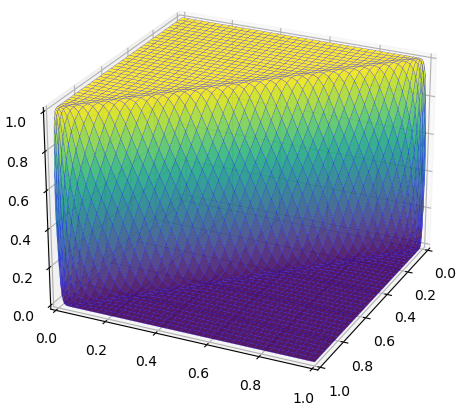}
    \end{minipage}
    \caption{Compares the numerical solution at $t=1.0$ for the $45\times 45$
    triangular and quadrilateral moving mesh.}
    \label{fig:mmesh-1}
\end{figure*}

\begin{figure*}[htbp]
    \centering
    \begin{minipage}{0.4\textwidth}
        \includegraphics[width=1.0\textwidth]{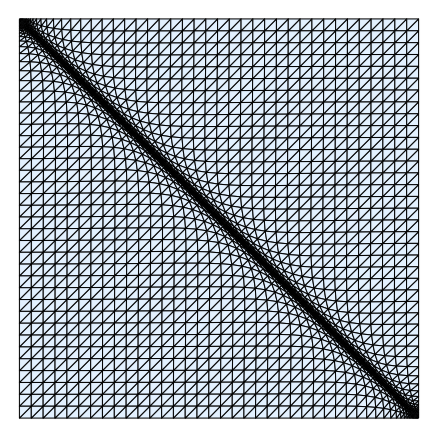}
    \end{minipage}
    \begin{minipage}{0.4\textwidth}
        \includegraphics[width=1.0\textwidth]{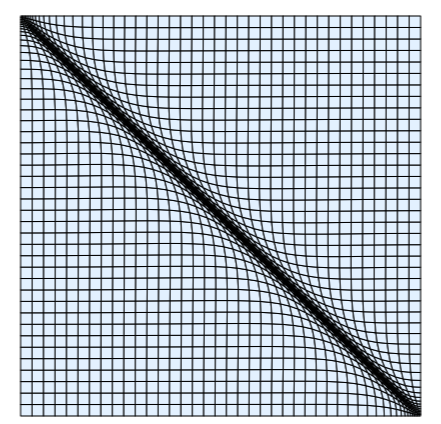}
    \end{minipage}
    \caption{Mesh distribution at $t=1.0$ for the $45\times 45$ triangular and
    quadrilateral moving mesh. Nodes are clearly concentrates within the sharp
    transition layer.}
    \label{fig:mmesh-2}
\end{figure*}

\paragraph{Benchmark.}
We apply the forward Euler time scheme with $a = 0.005$ and time step
$\Delta t = 0.001$.
The initial moving meshes include uniform triangular grids of sizes
$30\times30$, $45\times45$, and $60\times60$, as well as a $45\times45$
quadrilateral mesh.
For comparison, a fixed uniform triangular grid of size $160\times160$ is used.
See Table~\ref{tab:mmesh} for $L^2$ errors and times, Figure~\ref{fig:mmesh-1} and
\ref{fig:mmesh-2} for the solution and moving mesh respectively.
\begin{table}[htbp]
\centering
\begin{tabular}{clllllr}
\hline
\textbf{Type and Density}& \multicolumn{1}{c}{$t=0.5$} & \multicolumn{1}{c}{$t=0.75$} & \multicolumn{1}{c}{$t=1.0$} & \multicolumn{1}{c}{$t=1.25$} & \multicolumn{1}{c}{$t=1.5$} & \multicolumn{1}{c}{\textbf{Time (s)}} \\ \hline
Fix Tri $160\times 160$  & 0.0032488                   & 0.0040798                    & 0.0047631                   & 0.0041900                    & 0.0033862                   & 709.588                               \\
Tri $30\times 30$        & 0.0049148                   & 0.0051815                    & 0.0043073                   & 0.0046847                    & 0.0042645                   & 114.414                               \\
Tri $45\times 45$        & 0.0027004                   & 0.0030660                    & 0.0030052                   & 0.0029330                    & 0.0025511                   & 276.438                               \\
Tri $60\times 60$        & 0.0022169                   & 0.0026554                    & 0.0028656                   & 0.0025694                    & 0.0020576                   & 409.575                               \\
Quad $45\times 45$       & 0.0023646                   & 0.0024993                    & 0.0022470                   & 0.0019312                    & 0.0016308                   & 406.787                               \\ \hline
\end{tabular}
\caption{Convergence results ($L^2$ errors) for the Burgers' equation with moving mesh.}
\label{tab:mmesh}
\end{table}

\paragraph{Multiple backends.}
In the second test, we compare computational performance between CPU and GPU
execution.
With $a=0.005$ and $\Delta t=0.001$ over the time interval $[0, 0.25]$, we
measure the total runtime on triangular meshes of sizes
$60\times 60$, $80\times80$, $120\times120$, and $160\times160$.
The results, summarized in Table~\ref{tab:mmesh-time}, show that the
PyTorch-GPU backend achieves the best performance.

\begin{table}[htbp]
\centering
\begin{tabular}{ccrrrr}
\hline
\textbf{Device} & \textbf{Backend} & \multicolumn{1}{c}{$60\times60$} & \multicolumn{1}{c}{$80\times80$} & \multicolumn{1}{c}{$120\times120$} & \multicolumn{1}{c}{$160\times160$} \\ \hline
CPU             & NumPy            & 36.071                           & 107.114                          & 180.658                            & 387.741                            \\
CPU             & PyTorch          & 75.276                           & 121.256                          & 237.256                            & 526.509                            \\
GPU             & PyTorch          & 36.602                           & 48.130                           & 76.784                             & 102.323                            \\ \hline
\end{tabular}
\caption{Total runtime (seconds) for the Burgers' equation with different size of
moving meshes based on NumPy(CPU) and PyTorch(CPU/GPU).}
\label{tab:mmesh-time}
\end{table}

\subsection{Inverse Problem with Learning Automated FEM}
Electrical Impedance Tomography (EIT) serves as a typical example of PDE-based
inverse problems.
It aims to reconstruct the internal conductivity distribution $\sigma$ of a
domain $\Omega$ from boundary voltage and current measurements.
The governing physical model is described by the elliptic equation:
$$\nabla \cdot (\sigma \nabla u) = 0 \quad \text{in} \ \Omega,$$
with appropriate boundary conditions.
Recovering $\sigma$ from limited boundary data is notoriously ill-posed,
requiring innovative approaches that combine mathematical regularization with
data-driven learning.

\begin{figure}[htbp]
    \centering
    \includegraphics[width=0.80\textwidth]{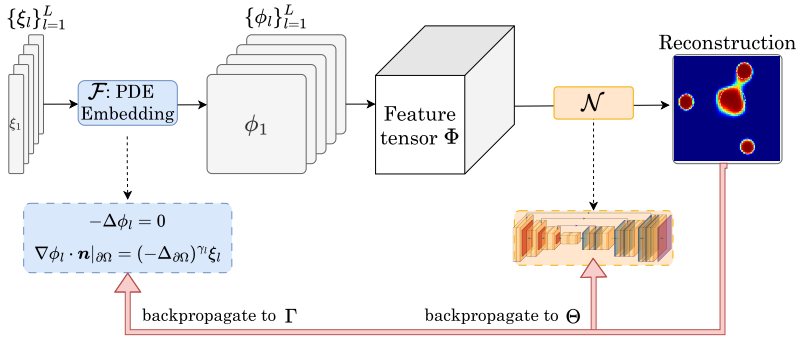}
    \caption{The $\gamma$-deepDSM framework.}
    \label{fig:gamma-deepDSM}
\end{figure}

\paragraph{The $\gamma$-deepDSM approach.}
Inspired by DDSM\cite{GUO2021}, a novel operator learning framework termed
$\gamma$-deepDSM\cite{ZHENG2024} was proposed
integrating a learnable preconditioning operator with a neural network for
enhanced reconstruction.
The method consists of two synergistic components (Figure~
\ref{fig:gamma-deepDSM}):
\begin{enumerate}
    \item Preprocessing via Learnable Fractional Laplace-Beltrami Operator
    \cite{MATEUSZ2017}.
    Boundary data $\xi_l = g_{D,l} - \Lambda_{\sigma_0} g_{N,l}$ is transformed
    using an fLB operator with learnable fractional orders $\gamma_l$:
    $$\mathfrak{L}_{\partial \Omega}^{\gamma_l} \xi_l = \sum_{k=1}^{K_0}
    \lambda_k^{\gamma_l} \alpha_k \psi_k,$$
    where $\lambda_k, \psi_k$ are eigenvalues and eigenfunctions of the
    Laplace-Beltrami operator.
    This step acts as a preconditioner to amplify high-frequency features
    critical for reconstruction.
    \item Feature Mapping via Neural Network.
    The preprocessed data $\{\phi_l\}_{l=1}^L$, obtained by solving:
    $$-\Delta \phi_l = 0 \quad \text{in} \ \Omega, \quad \nabla\phi_l\cdot\bm{n}
    = \mathfrak{L}_{\partial \Omega}^{\gamma_l} \xi_l \quad \text{on} \ \partial \Omega,$$
    are handled as a tensor and fed into a U-Net\cite{RONNEBERGER2015} to
    reconstruct the conductivity image $\sigma$.
\end{enumerate}

Thanks to the PyTorch support of FEALPy, we can easily implement the
$\gamma$-deepDSM framework by directly embedding FEM code into NN modules, as
shown in Listing~\ref{lst:gamma-deepdsm}.
\begin{lstlisting}[language=Python,
    caption={A PyTorch NN module containing a FEM solver. The 'self.solver'
    here is a wrapper for FEALPy code.},
    label={lst:gamma-deepdsm}
]
class DataFeature(nn.Module):
    ...
    def forward(self, input: Tensor) -> Tensor:
        BATCH, CHANNEL, NNBD = input.shape
        gnvn = gnvn.reshape(-1, NNBD) # [B*CH, NN_bd]
        val = self.solver.solve_from_current(gnvn)
        img = self.solver.value_on_nodes(val) # [B*CH, NN]
        return img.reshape(BATCH, CHANNEL, -1) # [B, CH, NN]
\end{lstlisting}

\paragraph{Results with FEALPy.}
Experimental results demonstrate that $\gamma$-deepDSM achieves approximately
15\% improvement in reconstruction accuracy compared to non-learnable approaches
across various noise levels and inclusion configurations.
The learned fractional orders $\gamma_l$ adapt to different frequency components
in the input data, with higher frequencies corresponding to smaller $\gamma_l$
values.

This example highlights FEALPy's unified tensor-based architecture that
seamlessly integrates traditional FEM with DL workflows.
The work creates a powerful framework for solving challenging inverse
problems while maintaining mathematical rigor and computational efficiency.

\subsection{Path Planning with Metaheuristic Optimization}
In addition to numerical calculations centered on solving direct or inverse PDE
problems, FEALPy also integrates a considerable number of intelligent
optimization algorithms in the \texttt{OPT} module.
Take path planning tasks as an example:
results of two problems, a 2D grid-based and a 3D terrain-based, are illustrated
in this subsection to demonstrate the adaptability of FEALPy for implementing
intelligent optimization algorithms based on tensor computing frameworks.

\begin{figure*}[htbp]
    \centering
    \begin{minipage}{0.45\textwidth}
        \includegraphics[width=1.0\textwidth]{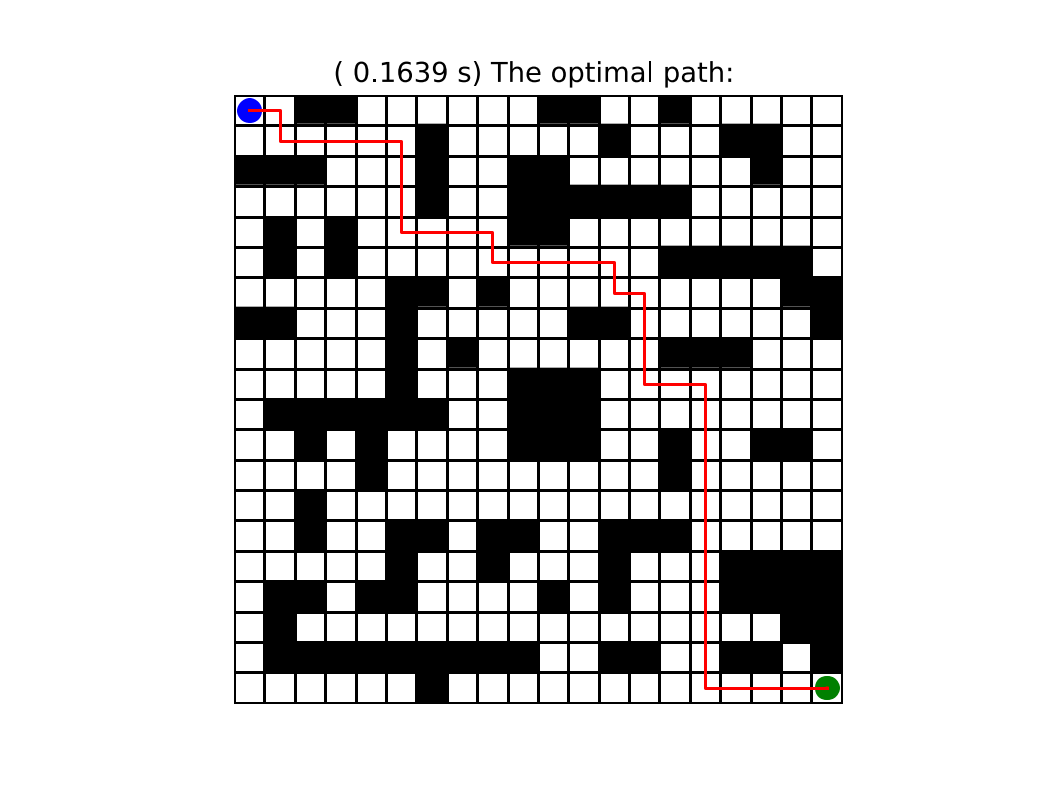}
        \caption{Path planning on a 2D $20\times 20$ grid map, where black
        squares denote obstacles and the red curve indicates the optimized
        route obtained by PSO.}
        \label{fig:path_planning}
    \end{minipage}
    \begin{minipage}{0.10\textwidth}
    \end{minipage}
    \begin{minipage}{0.45\textwidth}
        \includegraphics[width=1.0\textwidth]{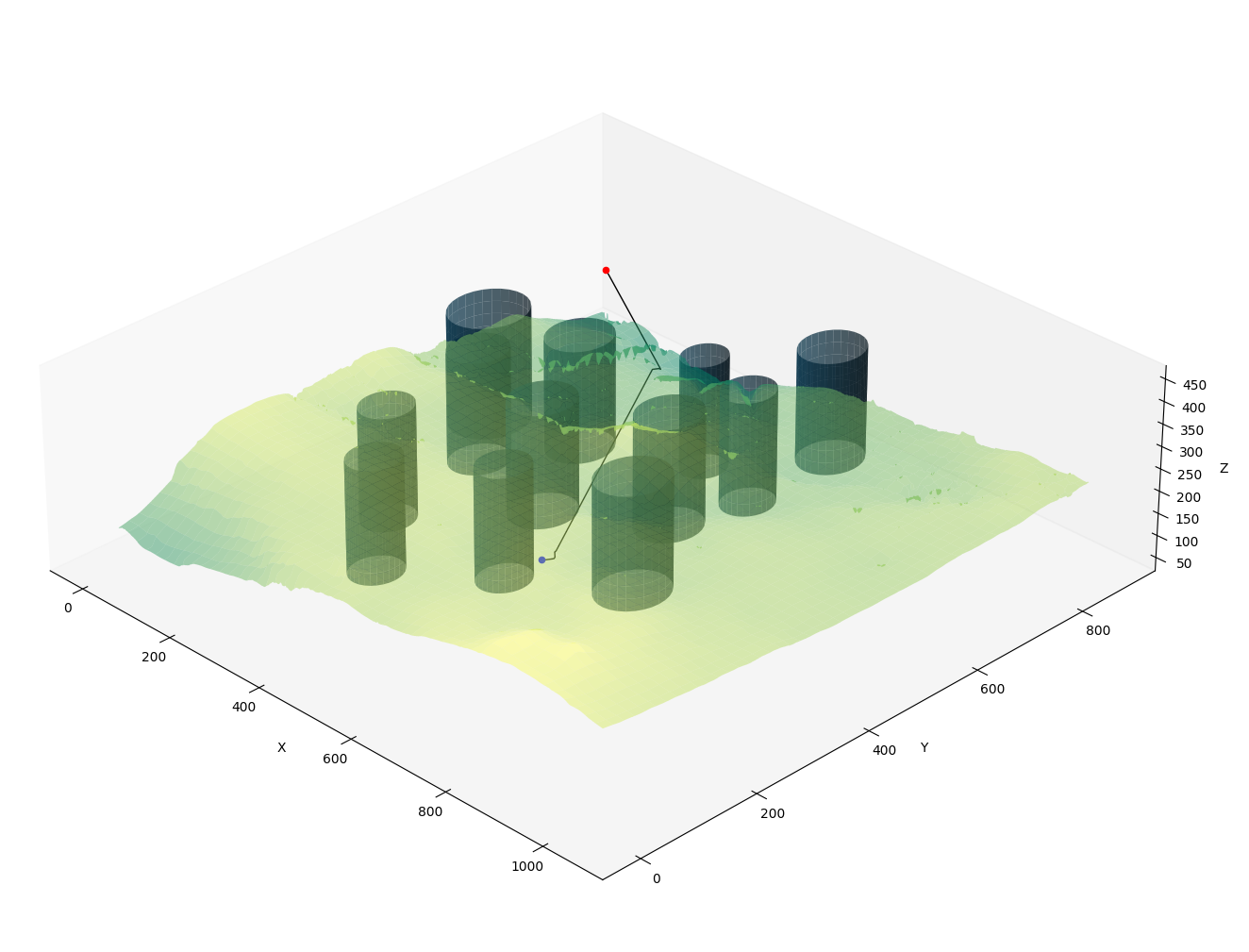}
        \caption{The optimized trajectory on a 3D terrain map, obtained by SAO.}
        \label{fig:path3d}
    \end{minipage}
\end{figure*}

\paragraph{PSO on a 2D grid map.}
A grid-based AGV (Automated Guided Vehicle) path planning problem is considered
on a $20\times 20$ map, where cells with value 1 represent obstacles and 0
represent free space.
Given the start point $(0,0)$ and end point $(19,19)$, we solve the problem
using the Particle Swarm Optimization (PSO)\cite{HAN2023} algorithm implemented
in FEALPy (Listing~\ref{lst:path2d}).
The resulting path is shown in Figure~\ref{fig:path_planning}.

\begin{lstlisting}[
    language=Python,
    caption={2D path planning model using PSO in FEALPy.OPT. Note that the
    PathPlanningModelManager generates conputational models
    (such as the planner in this code) as wrappers for
    the path planning scripts, which improved the readability and usability of
    the code in basic teaching. Readers can examine its internal implementation
    to learn how to use the ParticleSwarmOpt algorithm directly.},
    label={lst:path2d}
]
from fealpy.backend import bm
from fealpy.pathplanning.model import (
    PathPlanningModelManager,
)
from fealpy.opt import ParticleSwarmOpt

MAP = bm.array([
    [0, 0, 1, 1, 0, 0, ... , 0],
    [0, 0, 0, 0, 0, 0, ... , 0],
    ...
    [0, 0, 0, 0, 0, 0, ... , 0]
])
start_point, end_point = (0, 0), (19, 19)

# Configure model with PSO
options = {
    'MAP': MAP,
    'start_point': start_point,
    'end_point': end_point,
    'opt_alg': ParticleSwarmOpt
}
manager = PathPlanningModelManager('route_planning')
planner = manager.get_example(2, **options)
planner.solver()
planner.visualization()
\end{lstlisting}

\paragraph{SAO on a terrain map.}
A 3D terrain environment is constructed from the dataset
\texttt{ChrismasTerrain.tif}, with several cylindrical threat regions to be
avoided.
Given the start point $(800,100,150)$ and end point $(100,800,150)$,
the path was optimized using the Snow Ablation Optimization (SAO)\cite{DENG2023}
algorithm implemented in FEALPy.
Code for this problem is similar, but with a different map type (the thread
regions and terrain data), as shown in Listing~\ref{lst:path3d}.
The optimized trajectory in Figure~\ref{fig:path3d} shows the
3D terrain surface and feasible paths.

\begin{lstlisting}[
    language=Python,
    caption={3D terrain-based path planning using SAO in FEALPy.OPT.
    The parts that are the same as in the previous example have been omitted.},
    label={lst:path3d}
]
...
from fealpy.pathplanning import TerrainLoader
from fealpy.opt import SnowAblationOpt

terrain_data = TerrainLoader.load_terrain(
    'ChrismasTerrain.tif'
)
# [x, y, radius, height]
threats = bm.array([
    [400, 500, 200, 50], ..., [650, 450, 185, 50]
])
start_pos = bm.array([800, 100, 150])
end_pos   = bm.array([100, 800, 150])
options = {
    'threats': threats,
    'terrain_data': terrain_data,
    'start_pos': start_pos,
    'end_pos': end_pos,
    'opt_method': SnowAblationOpt
}
...
model = manager.get_example(1, **options)
solution, fitness = model.solver(n=10)
...
\end{lstlisting}

\section{Applications}

\subsection{SOPTX: High-Performance Multi-Backend Topology Optimization}
The SOPTX\cite{HE2025} framework was developed based on FEALPy primarily
targeting applications in structural topology optimization.
SOPTX fully inherits and extends FEALPy's multi-backends feature, diverse
numerical algorithm components, and general-purpose handling capabilities
for mesh and geometry.
Key technologies such as GPU parallel acceleration, sensitivity analysis
based on AD, and high-performance matrix assembly techniques are further
integrated based on this foundation, which significantly enhances the
computational efficiency of topology optimization problems while ensuring
flexibility and extensibility.

SOPTX is built upon a core architecture consisting of four main modules:
a Material Module, a Analysis Module, a Regularization Module, and an
Optimization Module.
These modules communicate and share data through well-defined interfaces,
forming a loosely coupled, easily extensible, and multi-backend framework for
topology optimization.
SOPTX has been successfully applied to solve various topology optimization
problems, generating results in excellent agreement with published literature
for classic benchmarks, such as the cantilever and MBB beams.
SOPTX is also demonstrated to be high performance, with a 67\% reduction in
assembly time via the high-performance matrix assembly technique, and a
8.1$\times$ overall speedup via GPU acceleration compared to CPU
in some large-scale 3D problems.

The modular architecture of SOPTX is poised for continued expansion:
\begin{itemize}
    \item Further enhance the boundary definition and geometric precision of
    optimization results by integrating technologies such as level-set and
    phase field solvers, along with adaptive mesh refinement.
    \item Extend from single-physics structural optimization to more complex,
    interdisciplinary applications, such as thermo-fluid-structure coupling,
    by combining with FEALPy's capabilities in solving multiphysics problems.
\end{itemize}

\subsection{FractureX: Flexible and Efficient Phase-Field Fracture Simulation}

As an application package built upon the FEALPy library, FractureX provides a
comprehensive and modular simulation environment for brittle fracture using the
phase-field method (PFM).
It is designed to bridge the gap between the theoretical sophistication of PFM
and its practical implementation, offering both computational robustness and
user flexibility for researchers in computational mechanics.

The core strength of FractureX lies in its modular architecture, which cleanly
decouples key physical and numerical components.
This design allows users to effortlessly combine different strain energy
decomposition models (e.g., isotropic, anisotropic, spectral, and hybrid
models \cite{MIEHE2010}),
energy degradation functions, and crack surface density functions
(AT1, AT2 \cite{WU2018}).
Such modularity facilitates rapid prototyping of novel constitutive models and
direct comparison between different modeling approaches within a unified
codebase.

FractureX adopts an automated adaptive mesh refinement (AMR) strategy driven by
a recovery-type error estimator\cite{TIAN2026}, which
significantly reduce the cost and improve the accuracy through dynamically
focusing the mesh resolution on the fracture region.
The main solver includes a robust Newton-Raphson scheme for the strongly
coupled displacement-phase-field system.
Inherited from FEALPy, GPU acceleration is ready for large-scale 3D simulations,
and AD is also available to skip conventional matrix assembly in rapid
tests for new models.

In summary, FractureX exemplifies the power of FEALPy as a foundation for
developing specialized scientific computing applications.
By abstracting the complexities of the FEM and providing a structured framework
for PFM, it lowers the barrier to entry for advanced fracture simulation and
serves as a versatile tool for both educational and state-of-the-art research
purposes.

\section{Discussion and Conclusions}

This work introduced FEALPy, a tensor-centric numerical simulation engine that
unifies traditional numerical simulation methods with modern tensor
computation/DL frameworks.
The open-source nature of FEALPy, along with its modular design, allows for
easier contributions from the community.

By organizing computations around a consistent tensor abstraction layer, FEALPy
achieves interoperability across multiple backends such as NumPy, PyTorch, and
JAX, allowing seamless execution on diverse hardware architectures including
CPUs and GPUs.
This design not only allows us researchers to leverage the latest advances in
hardware acceleration for tensors but also supports a wide range of software
ecosystems.

The tensor abstraction layer also results in the backend independence, which
allows developers to integrate future improvements.
Users and developers can make FEALPy adopt to new tensor libraries or user
C/C++ extensions for specific hardware platforms, by simply adding a new
backend adapter.
Newly-added tensor backends are available for all existing build-in/user
algorithms;
while new algorithms can run on all existing backends without modification
(Unless the algorithm takes advantage of a tensor library's unique features).

The modular architecture of FEALPy—comprising the tensor, commonality,
algorithmic, and field layers—effectively separates mathematical abstraction
from computational realization.
This design provides both clarity and flexibility for researchers to implement,
extend, and combine numerical algorithms without modifying the underlying data
structures.

Comprehensive numerical experiments have verified the correctness, efficiency,
extensibility and applicability, ranging from elasticity analysis and
high-order PDEs to moving mesh computation, inverse problems, and metaheuristic
optimization.
The derived application frameworks—SOPTX for topology optimization and
FractureX for phase-field fracture simulation—further illustrate how FEALPy
can serve as a foundation for developing multiple-backends, domain-specific
software.

All of these features together lower the migration barrier and learning cost
for researchers to adopt and contribute to FEALPy, fostering a collaborative
environment for advancing numerical simulation methods and their applications.
We have been expanding the collection of educational examples and
organizing the documentation of FEALPy to make it more accessible for students
and researchers new to the field, encouraging contributions from the community.

Future developments will focus on enhanced support for multi-physics coupling,
large-scale parallel computation, and data-driven numerical modeling.
With its unified tensor-based architecture and open, extensible design, FEALPy
aims to advance the convergence of numerical analysis, scientific computing,
and AI, providing a flexible platform for next-generation computational science.

\section*{Associated Content}

\subsection*{Data Availability Statement}

Source code of numerical examples are available at
\url{https://github.com/weihuayi/fealpy/tree/master/example}.
The user manual for FEALPy: \url{http://suanhai.com.cn/intro.html}.

\section*{Acknowledgments}
This work was supported by
the National Key R\&D Program of China (2024YFA1012600),
the National Natural Science Foundation of China (NSFC)
(Grant No. 12371410, 12261131501),
and the High Performance Computing Platform of Xiangtan University.

\bibliographystyle{plain}
\bibliography{FEALPyBib.bib}

@techreport{Chen:2008ifem,
  author = {Long Chen},
  title = {$i$FEM: an integrated finite element methods package in MATLAB},
  institution = {University of California, Irvine},
  year = {2009},
  url = {https://github.com/lyc102/ifem}
}

@article{Cai:2024afepack,
  author = {Z. Cai and R. Li and others},
  title = {AFEPack: A General-Purpose C++ Library for Numerical Solutions of
  Partial Differential Equations},
  journal = {Communications in Computational Physics},
  volume = {36},
  number = {1},
  pages = {274--318},
  year = {2024}
}

@article{Anderson:2021mfem,
  author = {R. Anderson and J. Andrej and A. Barker and J. Bramwell and
  J.-S.Camier and J. Cerveny and V. Dobrev and Y. Dudouit and A. Fisher and
  Tz. Kolev and W. Pazner and M. Stowell and V. Tomov and I. Akkerman and
  J. Dahm and D. Medina and S. Zampini},
  title = {MFEM: A Modular Finite Element Methods Library},
  journal = {Computers \& Mathematics with Applications},
  volume = {81},
  pages = {42--74},
  year = {2021},
  doi = {10.1016/j.camwa.2020.06.009}
}

@misc{BarattaEtal2023,
  title     = {{DOLFINx}: the next generation {FEniCS} problem solving environment},
  author    = {Baratta, Igor A. and Dean, Joseph P. and Dokken, J{\o}rgen S. and
  Habera, Michal and Hale, Jack S. and Richardson, Chris N. and
  Rognes, Marie E. and Scroggs, Matthew W. and Sime, Nathan and Wells, Garth N.},
  doi       = {10.5281/zenodo.10447666},
  year      = {2023},
  howpublished = {preprint}
}

@article{ScroggsEtal2022,
  title     = {Construction of arbitrary order finite element degree-of-freedom
  maps on polygonal and polyhedral cell meshes},
  author    = {Scroggs, Matthew W. and Dokken, J{\o}rgen S. and
  Richardson, Chris N. and Wells, Garth N.},
  journal   = {ACM Transactions on Mathematical Software},
  year      = {2022},
  volume    = {48},
  number    = {2},
  doi       = {10.1145/3524456},
  pages     = {{18:1--18:23}},
}

@article{BasixJoss,
  title     = {Basix: a runtime finite element basis evaluation library},
  author    = {Scroggs, Matthew W. and Baratta, Igor A. and
  Richardson, Chris N. and Wells, Garth N.},
  journal   = {Journal of Open Source Software},
  year      = {2022},
  volume    = {7},
  number    = {73},
  doi       = {10.21105/joss.03982},
  pages     = {3982}
}

@article{AlnaesEtal2014,
  title     = {Unified Form Language: A domain-specific language for weak
  formulations of partial differential equations},
  author    = {Alnaes, Martin S. and Logg, Anders and {\O}lgaard, Kristian B.
  and Rognes, Marie E. and Wells, Garth N.},
  journal   = {{ACM} Transactions on Mathematical Software},
  year      = {2014},
  volume    = {40},
  doi       = {10.1145/2566630},
}

@article{Arndt:2019dealII,
  author = {Daniel Arndt and Wolfgang Bangerth and Denis Davydov and
  Timo Heister and Luca Heltai and Martin Kronbichler and Matthias Maier and
  Jean-Paul Pelteret and Bruno Turcksin and David Wells},
  title = {The deal.II finite element library: design, features, and insights},
  journal = {Computers \& Mathematics with Applications},
  volume = {81},
  number = {3},
  pages = {407--422},
  year = {2021},
  doi = {10.1016/j.camwa.2020.02.022},
  eprint = {arXiv:1910.13247}
}

@article{Jasak:2009OpenFOAM,
  author = {H. Jasak and A. Jemcov and Z. Tukovic},
  title = {OpenFOAM: Open source CFD in research and industry},
  journal = {International Journal of Naval Architecture and Ocean Engineering},
  volume = {1},
  number = {2},
  pages = {89--94},
  year = {2009}
}

@article{         HARRIS2020-NUMPY,
    title         = {Array programming with {NumPy}},
    author        = {Charles R. Harris and K. Jarrod Millman and St{\'{e}}fan J.
                    van der Walt and Ralf Gommers and Pauli Virtanen and David
                    Cournapeau and Eric Wieser and Julian Taylor and Sebastian
                    Berg and Nathaniel J. Smith and Robert Kern and Matti Picus
                    and Stephan Hoyer and Marten H. van Kerkwijk and Matthew
                    Brett and Allan Haldane and Jaime Fern{\'{a}}ndez del
                    R{\'{i}}o and Mark Wiebe and Pearu Peterson and Pierre
                    G{\'{e}}rard-Marchant and Kevin Sheppard and Tyler Reddy and
                    Warren Weckesser and Hameer Abbasi and Christoph Gohlke and
                    Travis E. Oliphant},
    year          = {2020},
    month         = sep,
    journal       = {Nature},
    volume        = {585},
    number        = {7825},
    pages         = {357--362},
    doi           = {10.1038/s41586-020-2649-2},
    publisher     = {Springer Science and Business Media {LLC}},
    url           = {https://doi.org/10.1038/s41586-020-2649-2}
}

@article{PASZKE2019-PYTORCH,
    title={PyTorch: An Imperative Style, High-Performance Deep Learning Library},
    author={Adam Paszke and Sam Gross and Francisco Massa and Adam Lerer and James Bradbury and Gregory Chanan and Trevor Killeen and Zeming Lin and Natalia Gimelshein and Luca Antiga and Alban Desmaison and Andreas K{\"o}pf and Edward Yang and Zachary DeVito and Martin Raison and Alykhan Tejani and Sasank Chilamkurthy and Benoit Steiner and Lu Fang and Junjie Bai and Soumith Chintala},
    journal={ArXiv},
    year={2019},
    volume={abs/1912.01703},
    url={https://api.semanticscholar.org/CorpusID:202786778}
}

@software{JAX2018,
    author = {James Bradbury and Roy Frostig and Peter Hawkins and Matthew James Johnson and Chris Leary and Dougal Maclaurin and George Necula and Adam Paszke and Jake Vander{P}las and Skye Wanderman-{M}ilne and Qiao Zhang},
    title = {{JAX}: composable transformations of {P}ython+{N}um{P}y programs},
    url = {http://github.com/jax-ml/jax},
    version = {0.3.13},
    year = {2018},
}

@online{ARRAYAPI2024,
    title = {Python Array API Standard: https://data-apis.org/array-api/latest/},
    year = {2024},
    url = {https://data-apis.org/array-api/latest/},
    urldate = {2024-12-20}
}

@article{TENSORLY,
  author  = {Jean Kossaifi and Yannis Panagakis and Anima Anandkumar and Maja Pantic},
  title   = {TensorLy: Tensor Learning in Python},
  journal = {Journal of Machine Learning Research},
  year    = {2019},
  volume  = {20},
  number  = {26},
  pages   = {1-6},
  url     = {http://jmlr.org/papers/v20/18-277.html}
}

@article{CHEN2023GEO,
  title={Geometric Decomposition and Efficient Implementation of High Order Face and Edge Elements},
  author={Chen, Chunyu and Chen, Long and Huang, Xuehai and Wei, Huayi},
  journal={Communications in Computational Physics},
  volume={35},
  number={4},
  year={2024},
  publisher={Global-Science Press}
}

@inproceedings{RAVIART1977,
	address = {Berlin, Heidelberg},
	title = {A mixed finite element method for 2-nd order elliptic problems},
	isbn = {978-3-540-37158-8},
	booktitle = {Mathematical {Aspects} of {Finite} {Element} {Methods}},
	publisher = {Springer Berlin Heidelberg},
	author = {Raviart, P. A. and Thomas, J. M.},
	editor = {Galligani, Ilio and Magenes, Enrico},
	year = {1977},
	pages = {292--315},
}

@article{NEDELEC1986,
	title = {A new family of mixed finite elements in $R^3$},
	volume = {50},
	issn = {0945-3245},
	url = {https://doi.org/10.1007/BF01389668},
	doi = {10.1007/BF01389668},
	language = {en},
	number = {1},
	urldate = {2025-10-24},
	journal = {Numerische Mathematik},
	author = {Nédélec, J. C.},
	month = jan,
	year = {1986},
	pages = {57--81},
}

@article{BREZZI1985-BDM,
	title = {Two families of mixed finite elements for second order elliptic problems},
	volume = {47},
	issn = {0945-3245},
	url = {https://doi.org/10.1007/BF01389710},
	doi = {10.1007/BF01389710},
	language = {en},
	number = {2},
	urldate = {2025-10-24},
	journal = {Numerische Mathematik},
	author = {Brezzi, Franco and Douglas, Jim and Marini, L. D.},
	month = jun,
	year = {1985},
	pages = {217--235},
}

@article{HU2015-HuZhang,
	title = {A family of symmetric mixed finite elements for linear elasticity on tetrahedral grids},
	volume = {58},
	issn = {1869-1862},
	url = {https://doi.org/10.1007/s11425-014-4953-5},
	doi = {10.1007/s11425-014-4953-5},
	language = {en},
	number = {2},
	urldate = {2025-10-24},
	journal = {Science China Mathematics},
	author = {Hu, Jun and Zhang, ShangYou},
	month = feb,
	year = {2015},
	pages = {297--307},
}

@article{KIRBY2011-Bernstein,
	title = {Fast simplicial finite element algorithms using {Bernstein} polynomials},
	volume = {117},
	issn = {0945-3245},
	url = {https://doi.org/10.1007/s00211-010-0327-2},
	doi = {10.1007/s00211-010-0327-2},
	language = {en},
	number = {4},
	urldate = {2025-10-24},
	journal = {Numerische Mathematik},
	author = {Kirby, Robert C.},
	month = apr,
	year = {2011},
	pages = {631--652},
}

@article{WILLIAMS201418,
    title = {Symmetric quadrature rules for simplexes based on sphere close packed lattice arrangements},
    journal = {Journal of Computational and Applied Mathematics},
    volume = {266},
    pages = {18-38},
    year = {2014},
    issn = {0377-0427},
    author = {D.M. Williams and L. Shunn and A. Jameson},
}

@article{HUJUN2023,
    author = {Hu, Jun and Lin, Ting and Wu, Qingyu},
    year = {2023},
    month = {10},
    pages = {},
    title = {A Construction of $C^r$ Conforming Finite Element Spaces in Any Dimension},
    volume = {24},
    journal = {Foundations of Computational Mathematics},
    doi = {10.1007/s10208-023-09627-6}
}

@article{CHENLONG2021,
    author = {Chen, Long and Huang, Xuehai},
    year = {2021},
    month = {11},
    pages = {},
    title = {Geometric decompositions of the simplicial lattice and smooth finite elements in arbitrary dimension},
    doi = {10.48550/arXiv.2111.10712}
}

@article{CHEN2025,
  title={Implementation and Basis Construction for Smooth Finite Element Spaces},
  author={Chunyu Chen and Long Chen and Tingyi Gao and Xuehai Huang and Huayi Wei},
  journal={ArXiv},
  year={2025},
  volume={abs/2507.19732},
  url={https://api.semanticscholar.org/CorpusID:280323077}
}

@misc{ZHENG2024,
    title={$\gamma$-deepDSM for interface reconstruction: operator learning and a Learning-Automated FEM package},
    author={Yangyang Zheng and Huayi Wei and Shuhao Cao and Ruchi Guo},
    year={2024},
    eprint={2411.05341},
    archivePrefix={arXiv},
    primaryClass={math.NA},
    url={https://arxiv.org/abs/2411.05341},
}

@article{RUOLI2001,
    title = {Moving Mesh Methods in Multiple Dimensions Based on Harmonic Maps},
    journal = {Journal of Computational Physics},
    volume = {170},
    number = {2},
    pages = {562-588},
    year = {2001},
    issn = {0021-9991},
    doi = {https://doi.org/10.1006/jcph.2001.6749},
    url = {https://www.sciencedirect.com/science/article/pii/S002199910196749X},
    author = {Ruo Li and Tao Tang and Pingwen Zhang},
}

@article{RUOLI2002,
    title = {A Moving Mesh Finite Element Algorithm for Singular Problems in Two and Three Space Dimensions},
    journal = {Journal of Computational Physics},
    volume = {177},
    number = {2},
    pages = {365-393},
    year = {2002},
    issn = {0021-9991},
    doi = {https://doi.org/10.1006/jcph.2002.7002},
    url = {https://www.sciencedirect.com/science/article/pii/S0021999102970026},
    author = {Ruo Li and Tao Tang and Pingwen Zhang},
}

@article{HUANG2001,
    title = {Practical Aspects of Formulation and Solution of Moving Mesh Partial Differential Equations},
    journal = {Journal of Computational Physics},
    volume = {171},
    number = {2},
    pages = {753-775},
    year = {2001},
    issn = {0021-9991},
    doi = {https://doi.org/10.1006/jcph.2001.6809},
    url = {https://www.sciencedirect.com/science/article/pii/S0021999101968093},
    author = {Weizhang Huang},
}

@article{GUO2021,
    author = {Guo, Ruchi and Jiang, Jiahua},
    title = {Construct Deep Neural Networks based on Direct Sampling Methods for Solving Electrical Impedance Tomography},
    journal = {SIAM Journal on Scientific Computing},
    volume = {43},
    number = {3},
    pages = {B678-B711},
    year = {2021},
    doi = {10.1137/20M1367350},
    URL = {https://doi.org/10.1137/20M1367350},
    eprint = {https://doi.org/10.1137/20M1367350},
}

@article{MATEUSZ2017,
	title = {Ten Equivalent Definitions of the Fractional Laplace Operator},
	author = {Mateusz Kwaśnicki},
	journal = {Fractional Calculus and Applied Analysis},
	volume = {20},
	number = {1},
	pages = {7-51},
	year = {2017},
	issn = {1311-0454, 1314-2224}
}

@InProceedings{RONNEBERGER2015,
	author="Ronneberger, Olaf
	and Fischer, Philipp
	and Brox, Thomas",
	editor="Navab, Nassir
	and Hornegger, Joachim
	and Wells, William M.
	and Frangi, Alejandro F.",
	title="U-Net: Convolutional Networks for Biomedical Image Segmentation[J]",
	booktitle="Medical Image Computing and Computer-Assisted Intervention -- MICCAI 2015",
	year="2015",
	publisher="Springer International Publishing",
	address="Cham",
	pages="234-241"
}

@article{HAN2023,
    title = {A novel hybrid particle swarm optimization with marine predators},
    journal = {Swarm and Evolutionary Computation},
    volume = {83},
    pages = {101375},
    year = {2023},
    issn = {2210-6502},
    doi = {https://doi.org/10.1016/j.swevo.2023.101375},
    url = {https://www.sciencedirect.com/science/article/pii/S2210650223001487},
    author = {Baole Han and Baosheng Li and Chuandong Qin},
}

@article{DENG2023,
    title = {Snow ablation optimizer: A novel metaheuristic technique for numerical optimization and engineering design},
    journal = {Expert Systems with Applications},
    volume = {225},
    pages = {120069},
    year = {2023},
    issn = {0957-4174},
    doi = {https://doi.org/10.1016/j.eswa.2023.120069},
    url = {https://www.sciencedirect.com/science/article/pii/S0957417423005717},
    author = {Lingyun Deng and Sanyang Liu},
}

@misc{HE2025,
    title={SOPTX: A High-Performance Multi-Backend Framework for Topology Optimization}, 
    author={Liang He and Huayi Wei and Tian Tian},
    year={2025},
    eprint={2505.02438},
    archivePrefix={arXiv},
    primaryClass={math.AP},
    url={https://arxiv.org/abs/2505.02438},
}

@article{MIEHE2010,
    author = {Miehe, C. and Hofacker, M. and Welschinger, F.},
    title = {Phase field modeling of fracture in elastic-plastic solids: Application to ductile fracture},
    journal = {Computer Methods in Applied Mechanics and Engineering},
    year = {2010},
    volume = {199},
    pages = {2765--2778}
}

@article{WU2018,
    title = {A length scale insensitive phase-field damage model for brittle fracture},
    journal = {Journal of the Mechanics and Physics of Solids},
    volume = {119},
    pages = {20-42},
    year = {2018},
    issn = {0022-5096},
    doi = {https://doi.org/10.1016/j.jmps.2018.06.006},
    url = {https://www.sciencedirect.com/science/article/pii/S0022509618302643},
    author = {Jian-Ying Wu and Vinh Phu Nguyen},
}

@article{TIAN2026,
    title = {Adaptive finite element method for phase field fracture models based on recovery error estimates},
    journal = {Journal of Computational and Applied Mathematics},
    volume = {472},
    pages = {116732},
    year = {2026},
    issn = {0377-0427},
    doi = {https://doi.org/10.1016/j.cam.2025.116732},
    url = {https://www.sciencedirect.com/science/article/pii/S0377042725002468},
    author = {Tian Tian and Chunyu Chen and Liang He and Huayi Wei},
}

\end{document}